\documentclass[11pt]{amsart}

\usepackage[T1]{fontenc}
\usepackage[utf8]{inputenc}
\usepackage{lmodern}
\usepackage{stmaryrd}

\usepackage{amssymb}

\newtheorem{theorem}{Theorem}[section]
\newtheorem{lemma}[theorem]{Lemma}
\newtheorem{conjecture}[theorem]{Conjecture}
\newtheorem{corollary}[theorem]{Corollary}
\newtheorem{proposition}[theorem]{Proposition}

\theoremstyle{definition}
\newtheorem{definition}[theorem]{Definition}
\newtheorem{example}[theorem]{Example}

\theoremstyle{remark}
\newtheorem{remark}[theorem]{Remark}

\numberwithin{equation}{section}

\usepackage[american]{babel}
\usepackage[babel, final]{microtype}
\usepackage[all]{xy}
\usepackage{ifpdf}
\usepackage{comment}
\usepackage{multirow}
\usepackage{pinlabel}
\usepackage{geometry}
\usepackage[pdftex, dvipsnames]{xcolor}
\usepackage{amsmath,amsthm,amssymb,amsfonts,mathrsfs,tikz,tikz-cd,bm,verbatim,mathtools,caption,enumitem,csquotes,float}

\usepackage{amsmath,amsthm,amsfonts,amssymb,amscd,flafter,pinlabel}
\usepackage{mathtools}
\usepackage[mathscr]{eucal}
\usepackage{graphics}
 \usepackage[all]{xy}
\usepackage{caption, subcaption}
\usepackage{pgf}
\usepackage{tikz}
\usetikzlibrary{arrows,automata}

\usepackage[colorlinks=true, urlcolor=NavyBlue, linkcolor=NavyBlue, citecolor=NavyBlue, pdftitle={Constructions of symplectic forms on 4-manifolds}, pdfauthor={Peter Lambert-Cole}, pdfsubject={symplectic}, pdfkeywords={4-manifolds}]{hyperref}

\hypersetup{
    colorlinks=true,
    linkcolor=NavyBlue,
    citecolor=NavyBlue
}

\newcommand{\cR}{\mathcal{R}}
\newcommand{\cH}{\mathcal{H}}
\newcommand{\cA}{\mathcal{A}}

\newcommand{\ZZ}{\mathbb{Z}}
\newcommand{\QQ}{\mathbb{Q}}

\newcommand{\dbar}{\overline{\del}}

\newcommand{\betat}{\widetilde{\beta}}

\newcommand{\Ht}{\widetilde{H}}
\newcommand{\Mt}{\widetilde{M}}
\newcommand{\Rt}{\widetilde{R}}
\newcommand{\Sigmat}{\widetilde{\Sigma}}

\newcommand{\cT}{\mathcal{T}}
\newcommand{\cF}{\mathcal{F}}
\newcommand{\cK}{\mathcal{K}}
\newcommand{\cC}{\mathcal{C}}
\newcommand{\cD}{\mathcal{D}}

\newcommand{\cS}{\mathcal{S}}
\newcommand{\cB}{\mathcal{B}}

\newcommand{\CC}{\mathbb{C}}
\newcommand{\RR}{\mathbb{R}}

\newcommand{\CP}{\mathbb{CP}}

\newcommand{\del}{\partial}

\begin{document}

\title{Constructions of symplectic forms on 4-manifolds}
\date{August 2026}
\author{Peter Lambert-Cole}
\address{Department of Mathematics, University of Georgia, Athens, GA 30602}
\email{\href{plc@uga.edu}{plc@uga.edu}}

\subjclass[2010]{}

\dedicatory{}

\begin{abstract}

Given a symplectic 4-manifold $(X,\omega)$ with rational symplectic form, Auroux constructed branched coverings to $(\CP^2,\omega_{FS})$.  By modifying a previous construction of Lambert-Cole--Meier--Starkston, we prove that the branch locus in $\CP^2$ can be assumed holomorphic in a neighborhood of the spine of the standard trisection of $\CP^2$.  Consequently, the symplectic 4-manifold $(X,\omega)$ admits a cohomologous symplectic form that is K\"ahler in a neighborhood of the 2-skeleton of $X$.  We define the Picard group of holomorphic line bundles over the holomorphic 2-skeleton.
\end{abstract}

\maketitle



\section{Introduction}

\subsection{Motivation}

The motivation for this work is the geography problem for symplectic 4-manifolds and the symplectic BMY inequality.  The Bogomolev-Miyaoka-Yau inequality for complex surfaces states that if $X$ is a surface of general type, then
\begin{equation}
    \label{eq:BMY}
    c_1^2 \leq 3 c_2
\end{equation}
Miyaoka's proof proceeds via complex algebraic geometry \cite{Miyaoka}, while Yau's proof derives from his solution to the Calabi conjecture \cite{Yau}.  It is conjectured that the same inequality holds for closed symplectic 4-manifolds.  

\begin{conjecture}
    Let $(X,\omega)$ be a closed, minimal symplectic 4-manifold satisfying $c_1(X)^2  \geq 0$.  Then $c_1(X)^2 \leq 3c_2(X)$.
\end{conjecture}

One aim of this paper is to investigate whether it is possible to mimic these methods from complex algebraic geometry in the symplectic case.  The approach we take it to construct explicit symplectic forms that appear well-adapted to the smooth topological problem.  We give two constructions, of a so-called {\it K\"ahler trisection} and a {\it symplectic 1-cocycle}.  The underlying approach is the same in both cases, but produces slightly different geometric structures.

\subsection{Construction}

The constructions in this paper are slight variations of the construction in \cite{LMS}, with modifications according to the geometric structures of interest.

Given an integral symplectic form $\omega$ on $X$, Auroux constructed branched coverings $\pi_k: X \rightarrow \CP^2$ for $k \gg 0$ \cite{Auroux-branched}.  The symplectic form $k[\omega]$ can be recovered by pulling back the Fubini-Study form $\omega_{FS}$ by $\pi_k$ and perturbing to be nondegenerate.  The smooth topology of the branched covering depends {\it a priori} on the integral cohomology class of $\omega$ and is independent of choices for $k$ sufficiently large.  Auroux and Katzarkov strengthened this result, so that the branch locus in $\CP^2$ is a so-called quasiholomorphic curve \cite{AK-branched}.  Such a branch locus can be algebraically encoded by a factorization of the full twist in the braid group.  In \cite{LC-symp-surfaces}, it was shown that the branch locus $\cR$ in $\CP^2$ can be smoothly isotoped into transverse bridge position with respect to the standard trisection of $\CP^2$.  The trisection of $\CP^2$ then pulls back to a trisection of $X$, along with the relevant geometric structure.  

\subsection{Holomorphic trisections and K\"ahler-like forms}

An {\it almost-K\"ahler manifold} is a triple $(X,\omega,J)$ where $\omega$ is a symplectic form and $J$ is an almost-complex structure satisfying the compatibility condition $\omega(JX,JY) = \omega(X,Y)$ everywhere.  The pair $(\omega,J)$ determine a Riemannian metric $g(X,Y) = \omega(X,JY)$ on $X$.  The pair $(\omega,J)$ determines a K\"ahler structure if $J$ is integrable.

Our first main result is a construction of almost-K\"ahler structures that are integrable over the interesting part of the smooth topology of $X$.

\begin{corollary}
\label{cor:2-skeleton}
Let $(X,\omega)$ be a closed, symplectic 4-manifold with rational symplectic form.  There exists an almost-Kahler structure $(X,\omega',J)$ with $[\omega'] = [\omega] \in H^2_{DR}(X)$ and a smooth handle decomposition of $X$ such that $J$ is integrable over the 2-skeleton.
\end{corollary}

This follows immediately from the next theorem, as a trisection determines an `inside-out' smooth handle decomposition contained in the spine of the trisection.

\begin{theorem}
\label{thrm:main-kahler}
    Let $(X,\omega)$ be a closed, symplectic 4-manifold with rational symplectic form.  There exists an almost-Kahler structure $(X,\omega',J)$ with $[\omega'] = [\omega] \in H^2_{DR}(X)$ and a trisection $\cT$ of $X$ such that $J$ is integrable in a neighborhood of the spine $\cS$ of $\cT$.

    Moreover, there exists a decomposition $X = \nu(\cS) \cup Z$ such that
    \begin{enumerate}
        \item $(\nu(\cS),\omega',J)$ is K\"ahler with strictly pseudoconcave boundary,
        \item $Z$ is a (disconnected) subcritical Weinstein domain.
    \end{enumerate}
\end{theorem}

\begin{figure}
    \centering
    \includegraphics[width=.4\textwidth]{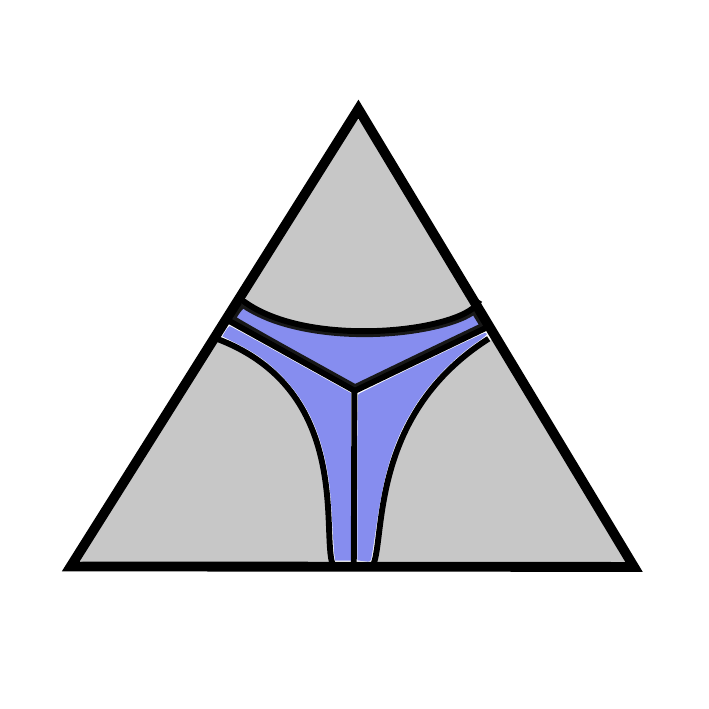}
    \caption{The decomposition of Theorem \ref{thrm:main-kahler} applied to the moment polytope of $\CP^2$} 
    \label{fig:CP2-decomp}
\end{figure}

There is a decomposition of $\CP^2$ as in Theorem \ref{thrm:main-kahler} depicted in terms of its moment polytope in Figure \ref{fig:CP2-decomp}.  The standard trisection of $\CP^2 = Z_1 \cup Z_2 \cup Z_3$ arises by taking the unit bidisk in each of the three standard affine charts.  Letting $Z'_{\lambda} \subset Z_{\lambda}$ be a small shrinking with strictly pseudoconvex boundary, we let $Z = Z'_1 \cup Z'_2 \cup Z'_3$ and $\nu(\cS)$ be its complement, which has strictly pseudoconcave boundary.  Pulling this decomposition back by a symplectic branched covering yields the decomposition of Theorem \ref{thrm:main-kahler}.

Since every subcritical Weinstein domain embeds symplectically in $(\CC^2,\omega_{std})$, we have the following corollary.

\begin{corollary}
\label{cor:main-kahler-cylinder}
    Let $(X,\omega)$ be a closed, symplectic 4-manifold with rational symplectic form.  There exists an almost-Kahler structure $(X,\omega',J)$ with $[\omega'] = [\omega] \in H^2_{DR}(X)$ such that $J$ is integrable, except on three trivial symplectic cylinders of the form $\left([a,b] \times Y_{\lambda}, d(e^t \alpha) \right)$, where $Y_{\lambda}$ is $S^3$ or $\#_k S^1 \times S^2$ for some $k > 0$.
\end{corollary}

\begin{remark}
    The subtlety that $[\omega']$ and $[\omega]$ are merely cohomologous is due to the fact that we appeal to the result of \cite{LC-symp-surfaces}, where the branch locus is only {\it smoothly} isotoped into transverse bridge position, as opposed to the Hamiltonian isotopy used in \cite{LMS}.
\end{remark}

\begin{remark}
    Viewed with respect to the inward-pointing coorientation, Theorem \ref{thrm:main-kahler} states that each boundary component admits a strongly pseudoconvex CR-structure.  If this CR-structure can be holomorphically embedded in $\CC^2$, then the K\"ahler structure could be extended across the subcritical Weinstein domain.  However, generic CR-structures on 3-manifolds do not holomorphically embed in $\CC^2$ and furthermore the existence of non-K\"ahler symplectic 4-manifolds implies that these boundary components cannot always be filled in holomorphically.
    \end{remark}

Since $J$ is integrable over the spine, we can attempt to mimic objects from algebraic geometry.  In particular, we introduce the Picard group $\text{Pic}(X,J,\cT)$ of holomorphic line bundles (Definition \ref{def:Picard}).  In fact, every isomorphism class of complex line bundle admits a holomorphic representative in some neighborhood of the 2-skeleton.

\begin{theorem}
\label{thrm:all-line-bundles}
Let $(X,\omega)$ be a closed symplectic 4-manifold with rational symplectic form.  There exists an almost-K\"ahler structure $(X,\omega',J)$, a trisection $\cT$ of $X$, and decomposition $X = \nu(\cS) \cup Z$ as in Theorem \ref{thrm:main-kahler} such that for any $\alpha \in H^2(X;\ZZ)$, there is some open subneighborhood $\nu(\cS)' \subset \nu(\cS)$ of the 2-skeleton and a holomorphic line bundle $L$ over $\nu(\cS)'$ such that $c_1(L) = \alpha$.
\end{theorem}

Neither Corollary \ref{cor:2-skeleton} nor Corollary \ref{cor:main-kahler-cylinder} require a smooth trisection decomposition, only Auroux's branched covering result.  One only needs to apply Lemma \ref{lemma:stein-approx} to the 1-skeleton of the branch curve to realize the holomorphic structures.  However, constructing a holomorphic structure on every smooth line bundle as in Theorem \ref{thrm:all-line-bundles} does require the smooth trisection.

We then attempt Hodge theory over the integrable region $\nu(\cS)$, viewed as a K\"ahler manifold with pseudoconcave boundary.  We can recover the Hodge-Riemann bilinear relations in subcritical dimensions (Theorem \ref{thrm:Hodge-Riemann}).

\subsection{Symplectic cone}

The {\it symplectic cone} of a smooth 4-manifold $X$ is the subset of $H^2_{DR}(X)$ that can be realized by symplectic forms.  Similarly, the {\it K\"ahler cone} of a complex surface $X$ is the subset of $H^{1,1}(X)$ that can be realized by K\"ahler forms.  Determining the symplectic cone of a general 4-manifold remains an open problem in general, although it is understood for some classes of 4-manifolds.  Seiberg-Witten invariants give some restrictions on the symplectic cone.  And determining the symplectic cone of $X \# k \overline{\CP}^2$ is connected to the problem of symplectic ball packing \cite{Biran-packing}.  In the K\"ahler case, the Nakai-Moishezon criterion for ampleness of holomorphic line bundles gives information about the K\"ahler cone.  There exists a symplectic Nakai-Moishezon-type criterion in the case when $b_2^+ = 1$ \cite{Biran-thesis,Biran-packing,Li-Liu}, but there is no corresponding general theorem.

Suppose that one of the holomorphic line bundles of Theorem \ref{thrm:all-line-bundles} admits many holomorphic sections over the spine $\nu(\cS)$.  The projective holomorphic embedding 
\[\iota: \nu(\cS) \hookrightarrow \CP^N\]
determined by this linear system gives a K\"ahler form $\iota^*(\omega_{FS})$ over the spine $\nu(\cS)$ that dominates the contact structures along the boundary $\del \nu(\cS)$.  In other words, it gives $\nu(\cS)$ the structure of a weak symplectic cap with $\omega$ in a prescribed cohomology class.  

If this form satisfies the stronger condition of being contact-type along the boundary --- in other words, being a {\it strong} symplectic cap --- then the boundary components can be filled in with symplectic 1-handlebodies to obtain a symplectic form $\omega$ on the closed 4-manifold $X$ with $\omega$ in a prescribed integral cohomology class of positive square.

However, it is not immediately obvious how to ensure this stronger condition or how to incorporate the restrictions from Seiberg-Witten theory.

\subsection{Symplectic 1-cocycles}

The notion of a {\it Weinstein trisection} of $(X,\omega)$ was introduced in \cite{LM-Rational}.  It is a trisection decomposition $X = X_1 \cup X_2 \cup X_3$ such that there exists a subcritical Weinstein structure on each sector $X_{\lambda}$ compatible with the global symplectic form.  Lambert-Cole, Meier and Starkston later showed that every closed, symplectic 4-manifold admits a Weinstein trisection \cite{LMS}.  A converse condition was introduced in \cite{LC-symp-cocycle} to ensure that a given trisection admits a global symplectic structure.

\begin{definition}
\label{def:symp-1-cocycle}
    Let $X$ be a closed, oriented 4-manifold and let $\cT$ be a trisection of $X$ with central surface $\Sigma$ and 3-dimensional sectors $H_1,H_2,H_3$.  A {\it symplectic 1-cocycle} for $(X,\cT)$ is a triple $(\betat_1,\betat_2,\betat_3)$ of closed, nonvanishing 1-forms on the three handlebodies such that the restrictions $\{\beta_{\lambda} = \betat_{\lambda|_{\Sigma}} \}$ satisfy the following conditions 
    \begin{enumerate}
        \item the cocycle condition:
        \[ \beta_1 + \beta_2 + \beta_3 = 0\]
        \item the pairwise positivity condition
        \[\beta_{\lambda} \wedge \beta_{\lambda+1} \geq 0\]
        \item they all vanish at the same $2g-2$ points $b_1,\dots,b_{2g-2}$ in $\Sigma$ and there exists a smooth function $f: \Sigma \rightarrow \RR$ and neighborhoods $\{b_i \in U_i\}$ such that $f|_{U_i} = \pm 1$ and $\sum_i f(b_i) = 0$ and
        \[\betat_{\lambda} = \beta_{\lambda} + d(sf)\]
        on a collar neighborhood of $\partial H_{\lambda} = \Sigma$.
    \end{enumerate}
\end{definition}

Recall that a boundary component of a symplectic manifold is {\it strongly concave} if there exists an inward-pointing Liouville vector field $\rho$.  This is equivalent to the criterion that the form $\alpha = \omega(\rho,-)$, which is a primitive for $\omega$, is a negative contact form on that boundary component (with respect to the {\it outward-pointing} coorientation).

\begin{theorem}
\label{thrm:main-symp-cocycle}
Let $X$ be a closed, oriented 4-manifold and let $\cT$ be a trisection of $X$ with central surface $\Sigma$.  Let $(\betat_1,\betat_2,\betat_3)$ be a symplectic 1-cocycle for $(X,\cT)$.  Then
\begin{enumerate}
    \item there is a symplectic form $\omega$ on $\nu(H_1 \cup H_2 \cup H_3)$ and three inward-pointing Liouville vector fields for $\omega$ along each of the three boundary components $\widehat{Y}_1, \widehat{Y}_2, \widehat{Y}_3$.
\item if each of the induced contact structures $(\widehat{Y}_{\lambda},\widehat{\xi}_{\lambda})$ are tight, then there exists a global symplectic form $\omega$ on $X$.
\end{enumerate}
\end{theorem}

In this paper, we prove two important existence statements about symplectic 1-cocycles and trisections.  In particular, we show that the criterion in \cite{LC-symp-cocycle} is completely general.

\begin{theorem}
\label{thrm:symp-trisection-existence-stabilization}
Let $(X,\omega)$ be a closed, symplectic 4-manifold.  Then
\begin{enumerate}
    \item there exists a trisection $\cT$ and symplectic 1-cocycle for $(X,\cT)$ such that all three induced contact structures on $\partial \nu(H_1 \cup H_2 \cup H_3)$ are tight.
    \item every stabilization of $\cT$ admits a symplectic 1-cocycle such that the induced contact structures on $\partial \nu(H_1 \cup H_2 \cup H_3)$ are tight.
\end{enumerate}
\end{theorem}

\subsection{Comparison}

A key fact is that any two smooth trisections $\cT_1,\cT_2$ of $X$ admit a common stabilization $\cT'$. The branch locus in the symplectic constructions described here can be perturbed to achieve a smooth stabilization of the trisection in the cover.  Consequently, we can directly compare two distinct symplectic structures on some $X$ via a common trisection.

\begin{theorem}
    Let $\omega_1,\omega_2$ be two rational symplectic forms on $X$.  There exists a single smooth trisection $\cT$ of $X$ such that the conclusions of Theorem \ref{thrm:main-kahler} and Theorem \ref{thrm:symp-trisection-existence-stabilization} hold for both $\omega_1$ and $\omega_2$, with respect to $\cT$.
\end{theorem}

\subsection{Acknowledgements}

I am deeply indebted to Phil Engel for explaining Miyaoka's proof to me in an intelligble way and to Cliff Taubes for pointing out a key error in an earlier draft.

\section{Trisections}
\label{sec:trisections}

Trisections are 4-dimensional generalization of Heegaard splittings of 3-manifolds.  In this section, we will describe the basic terminology and background results of the theory.  At this point, there is quite a voluminous literature defining the basic building blocks of the theory:
\begin{enumerate}
    \item trisections of 4-manifolds \cite{Gay-Kirby}
    \item relative trisections of 4-manifolds with boundary \cite{Relative-Tri}
    \item bridge trisections of surfaces \cite{MZ1,MZ2}
    \item constructing branched coverings via surfaces in bridge position \cite{LM-Rational,LMS}
    \item bridge trisections for properly embedded surfaces in 4-manifolds with boundary \cite{Meier-Relative-Bridge}
\end{enumerate}

The two key ideas we exploit in this paper are:
\begin{enumerate}
    \item Trisections pull back under branched coverings when the branch locus is in bridge position.  This has been used to construct trisections of projective surfaces \cite{LM-Rational} and symplectic 4-manifolds \cite{LMS}.
    \item In a simple branched covering, an (interior) perturbation of the branch locus downstairs results in an (interior) stabilization of the trisection upstairs.
\end{enumerate}

\subsection{Trisections}

Let $X$ be a closed, smooth 4-manifold.

\begin{definition}
A $(g;k_1,k_2,k_3)$-{\it trisection} $\cT$ of $X$ is a decomposition $X = Z_1 \cup Z_2 \cup Z_3$ such that for each $\lambda \in \ZZ_3$,
\begin{enumerate}
\item Each $Z_{\lambda}$ is diffeomorphic to $\natural_{k_{\lambda}} (S^1 \times B^3)$ for $k_{\lambda} \geq 0$,
\item Each double intersection $H_{\lambda} = Z_{\lambda-1} \cap Z_{\lambda}$ is a 3-dimensional, genus $g$ handlebody, and
\item the triple intersection $\Sigma = Z_1 \cap Z_2 \cap Z_3$ is a closed, genus $g$ surface.
\end{enumerate}
We refer to $\Sigma$ as the {\it central surface} of the trisection, the genus $g$ of $\Sigma$ as the {\it genus} of the trisection, and the union $H_1 \cup H_2 \cup H_3$ as the {\it spine} of the trisection.
\end{definition}

A trisection of $X$ can be encoded by a {\it trisection diagram} that describes how each $H_{\lambda}$ is attached to $\Sigma$ \cite[Definition 3]{Gay-Kirby}.  In particular, since $H_{\lambda}$ is a 1-handlebody, the gluing is determined by how a complete system of compressing disks is attached to $\Sigma$.  These compressing disks are attached along a cut system of curves in $\Sigma$ and these curves can be depicted on $\Sigma$ itself.

\subsection{Trisections of spheres}

The $n$-sphere has a standard trisection, obtained by pulling back a trisection on the unit disk $D$ in $\mathbb{R}^2$.  

\begin{definition}
The {\it standard trisection} $D = D^s_1 \cup D^s_2 \cup D^s_3 \subset \RR^2$ of the unit disk is defined by choosing the following subsets (in radial coordinates):
\begin{align*}
D^s_1 &= \left\{ 0 \leq \theta \leq \frac{2\pi}{3} \right\}, & D^s_2 &= \left\{\frac{2\pi}{3} \leq \theta \leq \frac{4\pi}{3} \right\}, &  D^s_3 &= \left\{ \frac{4\pi}{3} \leq \theta \leq 2 \pi \right\}.
\end{align*}
\end{definition}
It is symmetric under rotation by $2\pi/3$.

\begin{definition}[\cite{LC-Miller}]
\label{def:standard-tri-sphere}
The {\it standard trisection} of $S^{k}$ for $k\ge 2$ 
is the decomposition $\cT_{std} = \{\pi^{-1}(D^s_i) \cap S^{k}\}$ where $\pi: \mathbb{R}^{k+1} \rightarrow \mathbb{R}^2$ is a coordinate projection and $D^s_1 \cup D^s_2 \cup D^s_3$ is the standard trisection of the unit disk in $\mathbb{R}^2$.
\end{definition}

\begin{example}
\textit{The standard trisection of $S^2$}.  The triple intersection of $\cT_{std}$ consists of the North and South poles; the double intersections are three meridians connecting the poles; and the sectors are the three bigons bounded by the pairs of meridians.
\end{example}

\begin{example}
\label{ex:std-tri-s3}
\textit{The standard trisection of $S^3$}.  The triple intersection consists of a single $S^1$, unknotted in $S^3$; each double intersection is a $D^2$ Seifert disk for the unknot; and each sector of the trisection is a 3-ball.  In addition, we can give two other useful interpretations:
\begin{enumerate}
\item The lens space $L(1,3)$.  Here, $S^3$ is decomposed into three lenses: i.e take $D^2 \times [0,1]$ and collapse the vertical boundary to a circle.  These lenses are glued sequentially along their boundaries.
\item The genus 0 open book decomposition of $S^3$. Recall that $S^3$ has an open book decomposition $(B,\phi)$ where the binding $B$ is the unknot and the fibration $\phi: S^3 \smallsetminus \nu(B) \cong S^1 \times D^2 \rightarrow S^1$ is trivial.  This induces a trisection as follows.  Choose any three distinct pages $\left\{\phi^{-1}(\theta_1),\phi^{-1}(\theta_2), \phi^{-1}(\theta_3)\right\}$ of the open book decomposition.  Their union becomes the spine of the trisection and the binding $B$ is the codimension 2 stratum of the trisection.  The complement of the spine consists of three open 3-balls.
\end{enumerate}
\end{example}

The following lemma is immediate from the definition.

\begin{lemma}
\label{lemma:std}
Let $\cT_{std} = (S_1,S_2,S_3)$ be the standard trisection of $S^{k}$.  Then each $S_i$ is diffeomorphic to $B^{k}$; each double intersection $S_i \cap S_j$ is diffeomorphic to $B^{k-1}$ and the central surface is diffeomorphic to $S^{k-2}$.

In particular, the standard trisection of $S^4$ is a $(0;0)$-trisection.
\end{lemma}

\begin{example}
\label{ex:unbalanced-S4}
$S^4$ also admits an unbalanced $(1;1,0,0)$-trisection.  It can be described by a trisection diagram on $T^2$ with attaching curves $\alpha_1,\alpha_2$ parallel and the third attaching curve $\alpha_3$ geometrically dual to both.  This implies that $Y_1 = H_1 \cup -H_2$ is homeomorphic to $S^1 \times S^2$ and bounds $S^1 \times B^3$, while $Y_2 = H_2 \cup - H_3$ and $Y_3 = H_3 \cup -H_1$ are homeomorphic to $S^3$ and bound $B^4$.  By permuting the indices, we also obtain unbalanced $(1;0,1,0)$ and $(1;0,0,1)$-trisections of $S^4$.
\end{example}

\subsection{Connected sum and stabilization}

\begin{definition}
\label{def:connected-sum}
Let $X$ and $Y$ be closed, smooth, oriented 4-manifolds with trisections $\cT_X,\cT_Y$ given by the decompositions:
\[X = X_1 \cup X_2 \cup X_3 \qquad \text{and} \qquad Y = Y_1 \cup Y_2 \cup Y_3\]
The {\it connected sum} of $\cT_X$ and $\cT_Y$ is the trisection $\cT_Z$ of $Z = X \# Y$ obtained by performing the connected sum along points $x \in \Sigma_X \subset X$ and $y \in \Sigma_Y \subset Y$ in the central surfaces of the trisections.  Each sector is precisely $Z_{\lambda} =  X_{\lambda} \natural Y_{\lambda}$.  If $X$ has a $(g;k_1,k_2,k_3)$-trisection and $Y$ has a $(h;l_1,l_2,l_3)$-trisection, then this construction gives a $(g+h;k_1+l_1,k_2+l_2,k_3+l_3)$-trisection of $Z$.
\end{definition}

Recall from Example \ref{ex:unbalanced-S4} that $S^4$ admits an (unbalanced) trisection of genus 1, where the sector $X_{\lambda} \cong S^1 \times B^3$ and the remaining sectors are 4-balls.  Let $\cT_{\lambda}$ denote this trisection.

\begin{definition}
A {\it sector-$\lambda$ stabilization} of a trisection $\cT$ of $X$ is the trisection $\cT'$ obtained by taking the connected sum of $\cT$ with the trisection $\cT_{\lambda}$ of $S^4$
\end{definition}

Stabilization increases the trisection genus by 1.  Note that in \cite[Definition 8]{Gay-Kirby}, stabilizations are assumed to be {\it balanced} and increase the genus by 3;  however this operation is equivalent to three unbalanced stabilizations, one in each sector.

\begin{remark}
There is an equivalence between trisection decompositions and handle decompositions \cite[Section 4]{Gay-Kirby}.  Under this equivalence, an unbalanced stabilization of a trisection translates to adding a pair of canceling handles.
\end{remark}

\begin{theorem}[\cite{Gay-Kirby}] Let $X$ be a closed, smooth, oriented 4-manifold.  Then
\begin{enumerate}
\item $X$ admits a trisection,
\item any two trisections of $X$ become ambient isotopic after sufficiently many stabilizations, and
\item $X$ is determined up to diffeomorphism by the spine of any trisection $\cT$.
\end{enumerate}
\end{theorem}

\begin{remark}
All of the results in \cite{Gay-Kirby} assume {\it balanced} trisections, meaning $k_1 = k_2 = k_3$, although in practice it is often useful to consider unbalanced trisections.  It is always possible to balance a trisection by stabilization, so the distinction is irrelevant here.
\end{remark}

\subsection{Relative trisections}

Let $Z_0 = B^4$ and $Z_k = \natural_k S^1 \times B^3$ and let $Y_k = \del Z_k$.  For any integers $r,s$ let $\Sigma_{s,t}$ denote the compact, oriented surface of genus $s$ and $t$ boundary components.  For any integers $r,s,t$, let $H_{r,s,t}$ denote a {\it compression body} cobordism between $\Sigma_{r,t}$ and $\Sigma_{s,t}$.  In particular, if $r = s$, then $H_{r,s,t} = \Sigma_{r,t} \times [0,1]$ and if $r > s$, then $H_{r,s,t}$ is obtained from $H_{r,r,t}$ by attaching $r - s$ 3-dimensional 2-handles.  

Let $g,p,b$ be nonnegative integers with $g \geq p$ and $b \geq 1$.  Following \cite{Gay-Kirby,Relative-Tri}, there is a decomposition
\[Y_k = Y^-_{g,k;p,b} \cup Y^0_{g,k;p,b} \cup Y^+_{g,k;p,b}\]
where
\begin{align*}
Y^+_{g,k,p,b} & \cong H_{g,p,b} \\
Y^0_{g,k,p,b} & \cong \Sigma_{p,b} \times [0,1] \\
Y^-_{g,k,p,b} & \cong H_{g,p,b}
\end{align*}

\begin{definition}
Let $X$ be a compact, connected, oriented 4-manifold with connected boundary $Y$.  A $(g,k_1,k_2,k_3;p,b)$-{\it relative trisection} of $X$ is a decomposition
\[X = X_1 \cup X_2 \cup X_3\]
where
\begin{enumerate}
\item for each $\lambda = 1,2,3$, there is a diffeomorphism $\phi: X_{\lambda} \rightarrow Z_{k_{\lambda}} = \natural_{k_{\lambda}} S^1 \times B^3$ 
\item for each $\lambda = 1,2,3$ and taking indices mod 3, we have
\begin{align*}
\phi(X_{\lambda} \cap X_{\lambda - 1}) &= Y^-_{g,k_{\lambda},p,b} \\
\phi(X_{\lambda+1} \cap X_{\lambda}) &= Y^+_{g,k_{\lambda},p,b} \\
\phi(X_{\lambda} \cap \del X) &+ Y^0_{g,k_{\lambda},p,b}
\end{align*}
\end{enumerate}
Analogously to the closed case, we define $H_{\lambda} = X_{\lambda} \cap X_{\lambda-1}$, which is a compression body. 
\end{definition}

As in the closed case, we refer to $\Sigma = X_1 \cap X_2 \cap X_3$ as the {\it central surface}, the genus of $\Sigma$ as the {\it genus} of the trisection, and $H_! \cup H_2 \cup H_3$ as the {\it spine} of the trisection.
 
It follows from the construction (e.g. \cite[Lemma 11]{Relative-Tri}) that a relative trisection induces an open book decomposition of $\del X$.  Recall that an {\it open book decomposition} of 3-manifold $Y$ is a pair $(B,\pi)$ where $B$ is fibered link in $Y$ and $\pi: Y \setminus \nu(B) \rightarrow S^1$ is the fibration.

\begin{definition}
\label{def:std-tri-b4}
The {\it standard trisection} $\cT_{std}$ of $B^k$ for $k \geq 2$ is the decomposition $B^k = X_1 \cup X_2 \cup X_3$ where $X_{\lambda} = \pi^{-1}(D^s_{\lambda})$, where $\pi: \RR^k \rightarrow \RR^2$ is a coordinate projection, and $D = D^s_1 \cup D^s_2 \cup D^s_3$ is the standard trisection of the 2-disk.
\end{definition}

\begin{lemma}
Let $B^k = X_1 \cup X_2 \cup X_3$ be the standard trisection of $B^k$.  Then $X_{\lambda} \cap X_{\lambda+1}$ is diffeomorphic to $B^{k-1}$; the triple intersection $X_1 \cap X_2 \cap X_3$ is diffeomorphic to $B^{k-2}$; and $X_{\lambda} \cap \del B^k$ is diffeomorphic to $B^{k-1}$.

In particular, the standard trisection of $B^4$ is a $(0,0;0,1)$ relative trisection.
\end{lemma}

Note that if $X$ is a smooth 4-manifold with boundary with relative trisection $\cT_X$ and $Y$ is a closed, smooth 4-manifold with trisection $\cT_Y$, the {\it (interior) connected sum} of trisections $\cT_X \# \cT_Y$ of $X \# Y$ can be defined as in the closed case (Definition \ref{def:connected-sum}).

\begin{definition}
Let $X$ be smooth 4-manifold with boundary and let $\cT$ be a relative trisection of $X$. A {\it sector-$\lambda$ interior stabilization} of $\cT$ is the trisection obtained by taking the connected sum with the unbalanced, genus 1 trisection $\cT_{\lambda}$ of $S^4$.
\end{definition}

The main structural result for relative trisections is the following.

\begin{theorem}[\cite{Gay-Kirby}]
Let $X$ be a smooth, oriented 4-manifold with nonempty boundary.
\begin{enumerate}
\item For any fixed open book on the boundary $\del X$, there exists a relative trisection $\cT$  of $X$ inducing that open book deomposition.
\item Any two trisections of $X$ inducing the same open book on the boundary become ambient isotopic after a sequence of interior stabilizations.
\item the 4-manifold $X$ is determined up to diffeomorphism by the spine of a relative trisection.
\end{enumerate}
\end{theorem}

\subsection{Bridge position for surfaces}

Meier and Zupan generalized bridge position for links in a 3-manifold to bridge position for surfaces in a 4-manifold \cite{MZ1,MZ2}.  Meier has recently extended the theory of bridge trisections to properly embedded surfaces in 4-manifolds with boundary \cite{Meier-Relative-Bridge}.  We will give a brief overview, focusing only on the relevant case $X \cong B^4$.

A properly embedded tangle $\tau$ in 3-manifold with boundary $M$ is {\it trivial} if the arcs of $\tau$ can be simultaneously isotoped to lie in $\del M$.  A {\it disk-tangle} $\cD \subset X$ in a 4-manifold is a collection of properly embedded disks.  A disk-tangle is {\it trivial} if the disks can be simultaneously isotoped to lie in the boundary. All of the trivial tangles and disk-tangles we consider will be of the following form.  Let $\Sigma$ be a compact surface with nonempty boundary and let $Q$ be a collection of points in $\Sigma$.  Then
\begin{enumerate}
    \item $Q \times [0,1]$ is a trival tangle in $\Sigma \times [0,1]$, and
    \item $Q \times [0,1]^2$ is a trivial disk-tangle in $\Sigma \times [0,1]^2$
\end{enumerate}

Recall that a link $L \subset Y$ is in {\it bridge position} with respect to a Heegaard splitting $Y = H_1 \cup_{\Sigma} H_2$ if each intersection $\tau_{\lambda} = L \cap H_{\lambda}$ is a trivial tangle.  Extending, we say that a properly embedded, closed surface $\cK$ in a closed 4-manifold $X$ is in {\it bridge position} with respect to a trisection $\cT$ of $X$ if
\begin{enumerate}
    \item $\cK$ intersects the central surface $\Sigma$ of the trisection transversely in $2b$ points,
    \item $\cK$ intersects each handlebody $H_{\lambda}$ along a trivial tangle $\tau_{\lambda}$.
    \item $\cK$ intersects each 4-dimensional sector $X_{\lambda}$ along a trivial disk tangle.
\end{enumerate}
The parameter $b$ is the {\it bridge index} of the surface $\cK$.

Now suppose that $\cK$ is a properly embedded surface in $B^4$, let $\cT$ be a relative trisection of $B^4$, and suppose that $\del \cK$ is braided with respect to the induced open book decomposition on $\del B^4$.  

\begin{definition}
We say that $\cK$ is in {\it relative bridge position} with respect to $\cT$ if
\begin{enumerate}
    \item $\cK$ intersects the central surface $\Sigma$ transversely in $2b + n$ points,
    \item $\cK$ intersects each handlebody $H_{\lambda}$ along a trivial tangle, and
    \item $\cK$ intersects each sector $X_{\lambda}$ along a trivial disk-tangle.
\end{enumerate}
If $n$ is the braid index of $\del \cK$, we say that $b$ is the {\it bridge index} of $\cK$.
\end{definition}

\subsection{Bridge perturbation}
\label{sub:bridge-perturbation}

For motivation, recall that a {\it bridge perturbation} of a link $L \subset Y$ in bridge position is a local modification that increases the bridge index within the same ambient isotopy class.  We can describe the move abstractly and in coordinates:

\begin{enumerate}
    \item {\bf Abstractly} Near a bridge point $x \in \Sigma$, we can choose a local neighborhood $U$ and coordinates $(x,y,)$ on $U$ such that (a) the Heegaard surface intersects $U$ along the $xy$-plane, and (b) the link $L$ intersects $U$ along the $z$-axis.  Now choose an arc $\delta$ in $U$, such that one endpoint lies on $L$ and the other on $\Sigma$.  This determines a finger-move isotopy of $L$ through $\Sigma$.
    \item {\bf In coordinates}.  In the coordinates on $U$ from above, we can isotope $L$ to be the linear graph $L = (t,0,t)$.  By an isotopy, we can further isotope
 the link to be the cubic graph $L' = (t,0,t^3 - \epsilon t)$.
\end{enumerate}

We can define an analogous procedure for surfaces in bridge position.  Following \cite[Definition 9.9]{Meier-Relative-Bridge}, we first define a {\it finger-perturbation} abstractly.

\begin{definition}[\cite{Meier-Relative-Bridge}]
Let $\cK$ be a surface in relative bridge position with respect to a relative trisection $\cT$.  Let $x$ be a bridge point and $U$ a neighborhood of $X$ such that $\cK$ intersects $U$ along a small disk.  The intersection of $\cK$ with $U \cap H_{\lambda}$ is an arc $\tau_{\lambda}$ with one endpoint at $x$.  A {\it sector-$\lambda$ finger perturbation} of $\cK$ is the surface $\cK'$ obtained by a finger-move isotopy along a small arc $\delta$ in $U \cap H_{\lambda+1}$, with one endpoint on $\tau_{\lambda+1}$ and the other in $\Sigma$.
\end{definition}

The finger move increases the bridge index by 1; it increases the number of components of $\cK \cap X_{\lambda+1}$ by one (since the sector $X_{\lambda+1}$ of the trisection is opposite to the handlebody $H_{\lambda}$); and it fixes the number of components of $\cK \cap X_{\lambda-1}$ and $\cK \cap X_{\lambda}$.

We can also describe this move visually in 3 dimensions by projecting away one dimension.  Let $\cK$ be a surface in bridge position with respect to a trisection $\cT$.  We can locally choose coordinates $(x,y,z)$ on $\RR^3$ near a bridge point $\cK \cap \Sigma$ such that
\[ \cK = \{ z = 0\}\]
and the trisection is given by pulling back the standard trisection of the $xy$-plane by the projection $\pi: \RR^3 \rightarrow \RR^2$.

A bridge perturbation consists of locally replacing $\cK$ by
\[\cK' = \{ x = z^3 - yz \} \]
Intuitively, this is accomplished by {\it pleating} the surface $\cK$ through the center of the trisection.

Restricted to $\cK$, the projection $\pi$ is a homeomorphism.  Restricted to $\cK'$, however, there is a critical locus that projects to a semicubical cusp $(3t^2,-2t^3)$.  Over the interior of the cusp, the map is 3-to-1; over the exterior it is 1-to-1; at the origin there is a cusp singularity; and everywhere else is a fold singularity.

\subsection{Branched covers of surfaces in relative bridge position}

\begin{theorem}
\label{thrm:bridge-branch-tri}
Let $X$ be smooth 4-manifolds with boundary, let $\cK \subset X$ be a properly embedded surface, and let $\pi: \widetilde{X} \rightarrow X$ be an $n$-fold simple branched covering determined by a homomorphism $\rho: \pi_1(X \setminus \cK) \rightarrow S_n$.

\begin{enumerate}
    \item Suppose that $\cT$ is a relative trisection of $X$ and $\cK$ is in relative bridge position with respect to $\cT$.  Then $\widetilde{\cT} = \pi^{-1}(\cT)$ is a relative trisection of $\widetilde{X}$.
    \item Suppose that $\cK'$ isotopic to $\cK$ by a sector-$\lambda$ interior bridge perturbation.  Let $\widetilde{\cT},\widetilde{\cT}'$ be the induced relative trisections of $\widetilde{X}$.  Then $\widetilde{\cT}'$ is obtained from $\widetilde{\cT}$ by a sector-$(\lambda+1)$ interior stabilization.
    \end{enumerate}
\end{theorem}

\begin{proof}
Part (1) is the relative version of branched covering constructions of trisections, such as \cite[Theorem 4.1]{LMS}, \cite[Proposition 3.1]{LM-Rational}, and \cite[Proposition 13]{MZ2}.

An interior bridge perturbation increases the bridge index by 1 and the number of intersection points $\cK' \cap \Sigma$ by two.  Therefore, since the map $\pi$ is a simple branched covering, the Riemann-Hurwitz formula implies that the Euler characteristic of the central surface $\widetilde{\Sigma}$ upstairs decreases by 2 or equivalently the genus increases by 1.  The number of components of $\cK' \cap X_{\lambda-1}$ and $\cK' \cap X_{\lambda}$ stays the same, while the number of components of $\cK' \cap X_{\lambda+1}$ increases by 1.  Therefore, the Euler characteristic of $\widetilde{X}_{\lambda-1},\widetilde{X}_{\lambda}$ stay the same; since they are 1-handlebodies this implies that their smooth topology stays the same.  However, the Euler characteristic of $\widetilde{X}_{\lambda+1}$ decreases by 1.  Thus $\widetilde{X}'_{\lambda+1}$ is obtained from $\widetilde{X}_{\lambda+1}$ by attaching a 1-handle.  In summary, this exactly corresponds to the trisection stabilization.
\end{proof}

\subsection{Singular surfaces}

The branch loci for symplectic branched coverings $\pi: (X,\omega) \rightarrow (\CP^2,\omega_{FS})$ are singular surfaces, with self-transverse nodes and cusp singularities that are locally a cone on the right-handed trefoil.  We give a brief summary of \cite[Sections 3 and 4]{LMS}, which extends the (bridge) trisection machinery to singular surfaces.

A {\it singular knotted surface} $\cK$ in $X$ is the image of a smooth immersion, away from a finite number of critical points, such that all multiple points are transverse double points and in a small neighborhood of each critical value, the surface is diffeomorphic to the cone on a knot in $S^3$.  A {\it singular disk tangle} $\cD$ in $Z$ bounded by $L \subset Y = \partial D$ is a collection of properly embedded singular disks such that (a) each component of $\cD$ bounded by an unknot is smoothly embedded and boundary-parallel; (b) for  each component of $\cD$ with knotted boundary $K$, there exists an embedded 4-ball $B$ in $Z$ such that $B \cap Y$ is a 3-ball containing $K$ and the intersection of $\cD$ with $B$ is a smooth cone on $K$; (c) two components of $\cD$ are disjoint if their boundary knots are split and intersect transversely, away from the singular points, if their boundary knots are linked.  A singular surface $\cK$ is $X$ is in {\it bridge trisected position} with respect to a trisection if it is transverse to the central surface; the intersection $\tau_{\lambda} = H_{\lambda} \cap \cK$ is a boundary-parallel tangle; and the intersection of $\cK$ with each 4-dimensional sector is a singular disk tangle.

A {\it singular branched covering} $\pi: \widetilde{X} \rightarrow X$ with singular branch locus $\cK \subset X$ is a proper, smooth map such that (a) away from the singular points of $\cK$, the map is a branched covering, and (b) for each singular point $p$ of $\cK$, there exists a 4-ball neighborhood $U$ such that $\pi{^-1}(U)$ is a collection of $\widetilde{U}_! \cup \dots \cup \widetilde{U}_k$ and the restriction of $\pi$ to $\widetilde{U}_i$ is either a diffeomorphism or the cone on a branched covering $S^3 \rightarrow S^3$.  A trisection $\cT$ of $X$ pulls back to a trisection of $\widetilde{X}$ under a singular branched covering, provided that $\widetilde{X}$ is connected \cite[Theorem 4.1]{LMS}.

\section{Existence of symplectic 1-cocycles}

\subsection{Standard trisection of $(\CP^2,\omega_{FS})$}

The standard trisection of $\CP^2$ is induced by the toric structure (see citations).  In homogeneous coordinates $\{[z_1:z_2:z_3]\}$ on $\CP^2$, we have
\begin{align*}
    Z_{\lambda} &:= \{|z_{\lambda-1}| = 1 \, \text{ and } |z_{\lambda}|,|z_{\lambda+1}| \leq 1\} \\
    H_{\lambda} &:= \{|z_{\lambda-1}| = |z_{\lambda}| = 1 \, \text{ and } |z_{\lambda+1}| \leq 1\}\\
    \Sigma &= \{|z_1| = |z_2| = |z_3| = 1\}
\end{align*}
The central surface $\Sigma$ is a torus.  Since its Euler characteristic is 0, every harmonic 1-form is nonvanishing.  Therefore, symplectic data for this trisection consists of a triple of 1-forms $(\beta_1,\beta_2,\beta_3)$ on the handlebodies such that their restrictions to $\Sigma$ satisfy
\[\beta_{\lambda} \wedge \beta_{\lambda+1} > 0\]

These 1-forms can be obtained by tracing the Fubini-Study form $\omega_{FS}$ through the the Cech-DeRham isomorphism, as explained in \cite[Section 2.1]{LC-cohomology}.  With respect to polar coordinates $z_j = r_j e^{i \theta_j}$, this isomorphism sends $\omega_{FS}$ to the triple $(2 d \theta_1, 2 d \theta_2, 2(-d \theta_1 - \theta_2))$.

\subsection{Transverse bridge position and geometrically transverse surfaces}

\begin{definition}
\label{def:transverse-bridge}
Let $\cT$ be a trisection of $X$, let $\cK$ be a surface transverse to the trisection $\cT$, and let $(\beta_1,\beta_2.\beta_3)$ be a symplectic 1-cocycle for $\cT$.  \begin{enumerate}
    \item $\cK$ has {\it complex bridge points} if $\cK \cap \Sigma$ is disjoint from the zeros of $\beta$ and for each point of $\cK \pitchfork \Sigma$, there exists an integrable complex structure $J$ such that each $\betat_{\lambda}$ defines a $J$-holomorphic foliation on $H_{\lambda}$ and $\cK$ is locally $J$-holomorphic.
    \item $\cK$ is {\it geometrically transverse} if it has complex bridge points and $\betat_{\lambda}(\tau'_{\lambda}) > 0$ everywhere.
\end{enumerate}
If $\cK$ is in bridge position and geometrically transverse, we say it is in {\it transverse bridge position}.
\end{definition}
In other words, a surface is geometrically transverse if it has complex bridge points and each tangle $\tau_{\lambda}$ is positively transverse to the foliation of $H_{\lambda}$ defined by $\betat_{\lambda}$.

\begin{lemma}
Let $\cK$ be a symplectic surface (possibly nodal and/or cuspidal) that is geometrically transverse.  There exists a sector-$\lambda$ finger perturbation of $\cK$ to $\cK'$, supported in a neighborhood of some bridge point, such that $\cK'$ is also geometrically transverse.
\end{lemma}

\begin{proof}
Without loss of generality, we assume $\lambda = 1$ and the bridge point is a positive intersection between $\cK$ and $\Sigma$.  We can locally assume that $\beta_1 = dy$ and $\beta_2 = -dx$ and $\beta_3 = dx - dy$ and the tangles $\tau_1,\tau_2,\tau_3$ are incident to the bridge point $b = (0,0)$.  For some small $\epsilon > 0$, add a positive bridge point $b_+$ at $(-2\epsilon,-\epsilon)$ and a negative bridge point $b_-$ at $(-\epsilon,\epsilon)$.  Split $\tau_1$ into two oriented arcs, one ending at $b_+$ and the other connecting $b_-$ to $b$; we can assume these satisfy $\beta_1(\tau'_1) = dy(\tau'_1) > 0$.  Then we can add oriented arcs of $\tau_2$ and $\tau_3$ from $b_-$ to $b_+$, which satisfy the positivity condition.  Compare with Figure \ref{fig:mini-stabilization-crossing}.
\end{proof}

\subsection{Existence of branch locus}

Auroux, building on techniques introduced by Donaldson, proved that every closed symplectic 4-manifold admits a branched covering over $(\CP^2,\omega_{FS})$ \cite{Auroux-branched}.  This was strengthened by Auroux-Katzarkov to show that the branch curve in $\CP^2$ is `quasiholomorphic' \cite{AK-branched}.  Intuitively, this means the surface intersects a generic pencil on $\CP^2$ akin to a genuine holomorphic curve.  In \cite{LMS}, this essential geometric property was called a 'symplectic braided surface'.  However, precise definitions are irrelevant for the current application.  The important point is that the branch locus can be algebraically encoded, up to isotopy, by words in the braid group $B_n$, where $n$ is the degree of the branch locus in $\CP^2$.

\begin{theorem}[\cite{AK-branched}]
Every closed symplectic 4-manifold admits a symplectic singular branched covering over $\CP^2$, such that the branch locus in $\CP^2$ is a (singular) symplectic braided surface.
\end{theorem}

\begin{theorem}[\cite{LC-symp-surfaces}]
    Let $\cK$ be a singular, symplectic braided surface in $\CP^2$ with $A_k$-singularities.  There exists a smooth isotopy of $\cK$ such that it is in transverse bridge position with respect to the standard trisection of $\CP^2$.
\end{theorem}

\begin{remark}
Great care was taken in \cite{LMS} to ensure that the branch locus in $\CP^2$ was isotoped symplectically, which preserves the symplectic structure obtained by pulling back.  However, since we are only interested in the diffeomorphism type and will construct a potentially new symplectic form $\widetilde{\omega}$ anyway, we just need to ensure that the branch locus is smoothly isotopic to the given curve.
\end{remark}

Combining these two results, we have the following corollary.

\begin{corollary}
Let $(X,\omega)$ be a closed symplectic 4-manifold.  Then there exists a nodal, cuspidal surface $F$ in transverse bridge position with respect to the standard trisection of $\CP^2$ such that $X$ is the branched cover of $\CP^2$ over $F$.
\end{corollary}

For completeness, we review the construction here.  Surfaces in $\CP^2$ can be encoded by {\it torus diagrams} on the central surface $\Sigma$ of the standard trisection.  If $\cK$ is in general position with respect to the trisection, it intersects $\Sigma$ in a finite number of points and each solid torus $H_{\lambda}$ along a tangle $\tau_{\lambda}$.  Generically, this tangle misses the core $B_{\lambda}$ of $H_{\lambda}$ and there is a projection $H_{\lambda} \smallsetminus B_{\lambda} \cong \Sigma \times (0,1] \rightarrow \Sigma$.  The tangles $\tau_1,\tau_2,\tau_3$ project onto collections of arcs that we label $\cA,\cB,\cC$, respectively.  See \cite[Section 2.2]{LC-Thom}.

Let $\Sigma = T^2$ be the central surface of the trisection of $\CP^2$.  We have coordinates $(x,y)$ for $x,y \in [0,1]$ such that the foliations are 
\[\beta_1 = dy \qquad \beta_2 = -dx \qquad \beta_3 = dx - dy\]
We can extend these foliations across the handlebodies.  Let $\alpha_1 = (x,0)$ and $\alpha_2 = (0,x)$ and $\alpha_3 = (-x,-x)$ for $x \in [0,1]$; these curves bound disks in $H_1$ and $H_2$ and $H_3$, respectively, that are leaves of the foliations $\cF_1$ and $\cF_2$ and $\cF_3$, respectively. Let $p = (0,0) = \alpha_1 \cap \alpha_2$ be a fixed basepoint.

A key result of \cite{AK-branched} is that the branch locus $\cR$ is encoded algebraically by a braid factorization of the full twist.  Specifically, if $\cR$ is a degree $d$ surface in $\CP^2$, it is determined by a factorization
\[\Delta^2_d = (g_1 \sigma^i_1 g^{-1}_1) \cdots (g_n \sigma^{i_n}_1 g^{-1}_n)\]
where $\sigma_1$ is the first half twist in the Artin braid group and $i_j \in \{-2,1,2,3\}$.  Reconstructing $\cR$ from this factorization is described in \cite{AK-branched,LC-symp-surfaces}.  

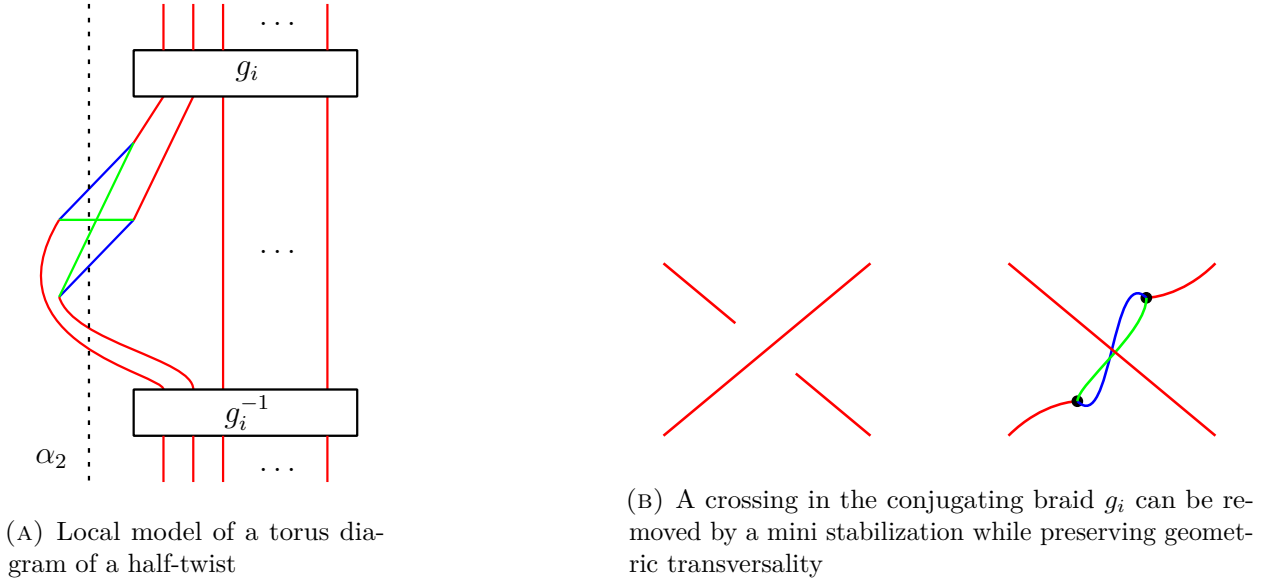
\begin{figure}[h!]
\centering
	\begin{subfigure}[b]{0.3\textwidth}
		\labellist
			\large\hair 2pt
			\pinlabel $g_i$ at 220 396
			\pinlabel $g^{-1}_i$ at 220 82
			\pinlabel $\dots$ at 250 440
			\pinlabel $\dots$ at 250 230
			\pinlabel $\dots$ at 250 30
			\pinlabel $\alpha_2$ at 40 40
		\endlabellist
		\includegraphics[width=\textwidth]{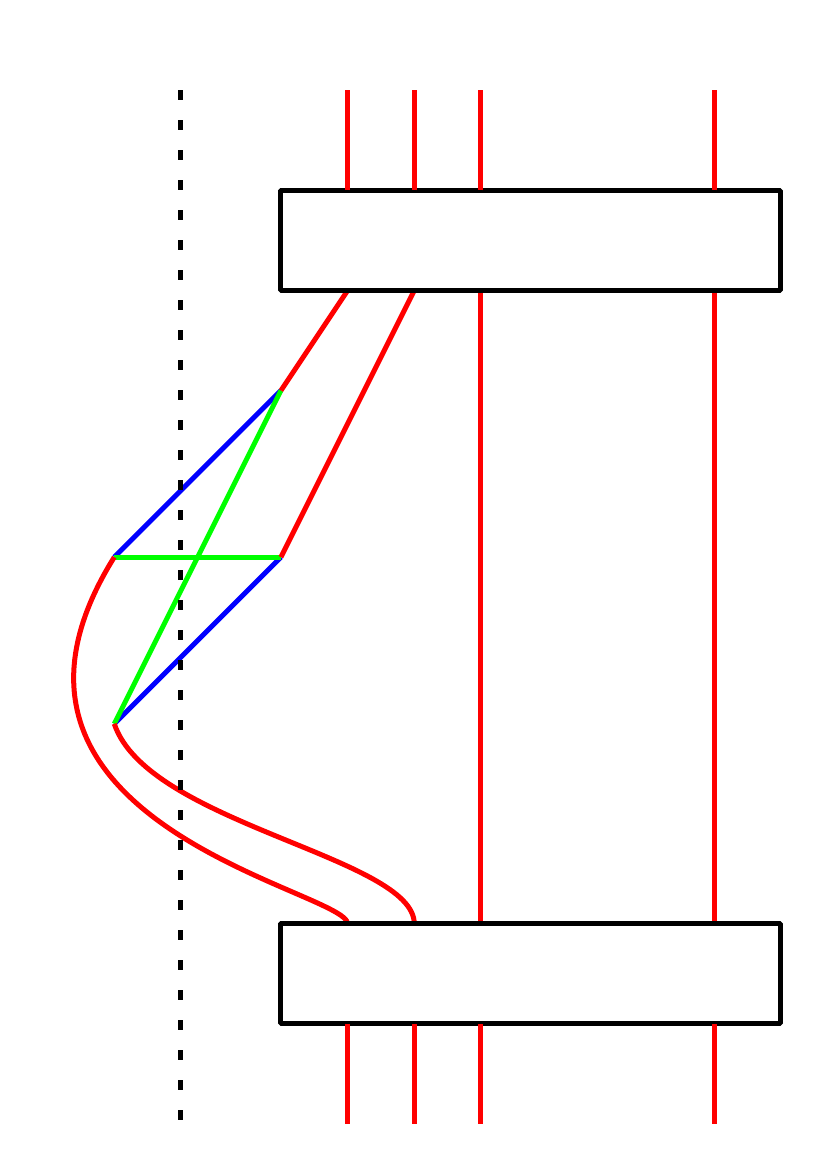}
		\caption{Local model of a torus diagram of a half-twist}
		\label{fig:band-local-model}
	\end{subfigure}
	\hfill
	\begin{subfigure}[b]{0.5\textwidth}
		\includegraphics[width=\textwidth]{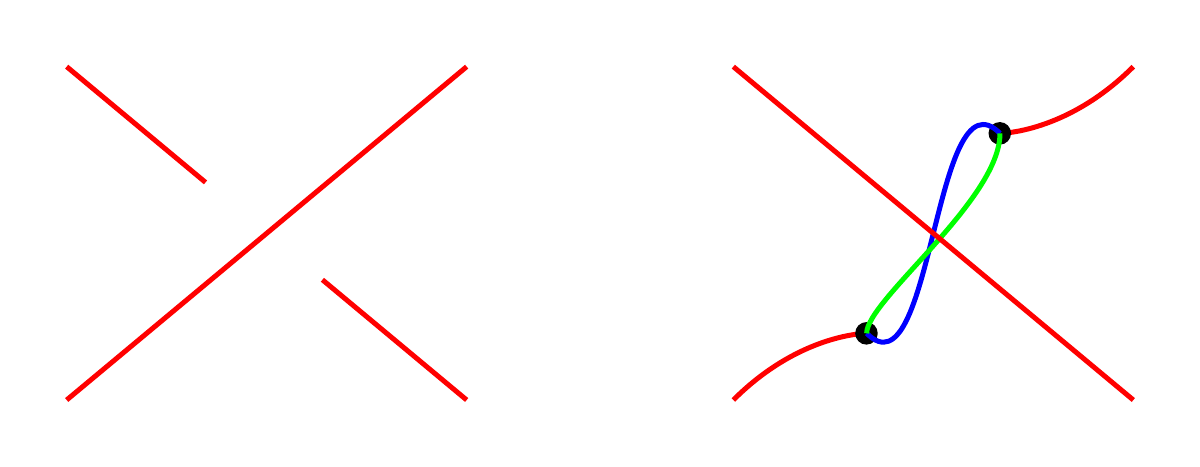}
		\caption{A crossing in the conjugating braid $g_i$ can be removed by a mini stabilization while preserving geometric transversality}
		\label{fig:mini-stabilization-crossing}
	\end{subfigure}
\caption{Component pieces of a torus diagram for the branch locus $\cR$}
\label{fig:branch-diagram}
\end{figure}

A torus diagram for $\cR$ in transverse bridge position in constructed in \cite[Proposition 4.6]{LC-symp-surfaces}.  It is determined by stacking local models for each conjugate half-twist $g_j \sigma^{i_j}_1 g^{-1}_j$, as  in Figure \ref{fig:band-local-model}.  Crossings in the conjugating braid $g_i,g^{-1}_i$ can be removed by a mini-stabilization (Figure \ref{fig:mini-stabilization-crossing}).  The only modification we require here is to dictate that each $\cB$ arc in the local model of a half-twist intersects the curve $\alpha_2$ transversely once.  In each sector, we obtain a transverse link $R_{\lambda} = \cR \cap Y_{\lambda}$.  Each component of $R_{\lambda}$ is either an (a) unknot, (b) Hopf link or (c) right-handed trefoil.  The surface $\cR$ is obtained by capping off each component of $R_{\lambda}$ with either (a) trivial disk, (b) pair of transversely intersecting trivial disks, or (c) a cone on a right-handed trefoil.

\subsection{Pulling back symplectic 1-cocycles}

A symplectic 1-cocycle $\{\betat_1,\betat_2,\betat_3)$ on a trisection $\cT$ of $X$ and a branch locus $\cK$ in singular, transverse bridge position pulls back under the branched covering
\[\pi: (\widehat{X},\widehat{\cT}) \rightarrow (X,\cT)\]

\begin{enumerate}
\item Away from the branch locus $\tau$, the map $\pi$ is a local diffeomorphism so $\pi^*(\betat_{\lambda})$ remains a closed, nonvanishing 1-form.  Similarly, along $\Sigma$ we can pull back the restrictions $\{\beta_{\lambda}\}$.
\item Along the branch locus $\tau_{\lambda}$, we have that $\betat_{\lambda}(\tau'_{\lambda}) > 0$ by assumption, since $\cK$ is in transverse bridge position.  We can therefore locally choose coordinates $(x,y,z)$ such that $\betat_{\lambda} = dz$ and $\tau = (0,0,t)$ for $t \in (-\epsilon,\epsilon)$.  The branched covering map is therefore given by 
\[\pi: (x,y,z) \mapsto (x^2 - y^2, 2xy, z)\]
and in particular $\pi^*(dz) = dz$ and the pullback of $\betat_{\lambda}$ is a closed, nonvanishing 1-form.
\item Finally, we check near the bridge points.  By assumption, we can choose a local neighborhood and identify it with $\CC^2$ with coordinates $(z_j = x_j + i y_j)$ such that
\begin{align*}
    \cK &= \{ z_2 = 0\} & \Sigma &= \{z_2 = \epsilon  \overline{z_1}\}
\end{align*}
where $\epsilon \in \{\pm 1\}$ is the sign of the intersection of $\cK$ with $\Sigma$.  Set 
\begin{align*}
    t &= x_1 - x_2 & s &= -y_1 - y_2 \\
    \betat_1 &= dy_1 - dy_2 & \betat_2 &= -dx_1 - dx_2
\end{align*}
Then the trisection decomposition is given by
\begin{align*}
H_1 &= \{s \geq 0, t = 0\} & Z_1 & =\text{max}(-s,t) \leq 0 \} \\
H_2 &= \{t \geq 0, -t -s = 0\} & Z_2 &=\text{max}(-t,-t-s) \leq 0\} \\
H_3 &= \{-t - s \geq 0, s = 0\} & Z_3 &= \text{max}\{(t+s,s) \leq 9\}
\end{align*}
and
\[dt \wedge \betat_1 - ds \wedge \betat_2 = 2 \left( dx_1 \wedge dy_1 + dx_2 \wedge dy_2 \right) = 2 \omega_{std}\]

The branched covering is locally given by the map $\pi: (z_1,z_2) \rightarrow (z_1,z_2^2)$.  The pulled-back trisection $\pi^*(\cT)$ is determined by the functions
\begin{align*}
t &= x_1 + \epsilon(- x_2^2 + y_2^2) & s &= -y_1 - \epsilon(2x_2 y_2)
\end{align*}
and the 1-forms pull back as
\begin{align*}
    \betat_1 &= dy_1 - 2\epsilon d(x_2y_2)\\
    \betat_2 &= -dx_1 + 2\epsilon d(y_2^2 - x_2^2)
\end{align*}
Restricting this forms to $\Sigma = \{t = s = 0\}$, we have that
\begin{align*}
   \beta_1 = \betat_1|_{\Sigma} &= -4\epsilon d(x_2y_2) \\
   \beta_2 = \betat_2|_{\Sigma} &= -4\epsilon d(y_2^2 - x_2^2)
\end{align*}
Using $x_2,y_2$ as coordinates along $\Sigma$ we see that
\[\beta_1 \wedge \beta_2 = 16\epsilon^2(x_2^2 + y_2^2) dx_2 \wedge dy_2 \geq 0\]
and
\[ \betat_1 = \beta_1 + d(s \cdot \epsilon) \qquad \qquad \betat_2 = \beta_2 + d(t \cdot \epsilon)\]
Consequently, the form $\beta_1$ has a hyperbolic zero at the branch locus. 
 Moreover, we can set $f = \epsilon \in \{\pm 1\}$ near the branch points in $\pi^*(\Sigma)$, according to the sign of the intersection $\cK \pitchfork \Sigma$.  Since $\Sigma$ is nullhomologous, as it bounds the handlebody $H_{\lambda}$, the algebraic count of intersections is 0.
\end{enumerate}

\subsection{Branched covers over transverse knots}

The final piece in proving Theorem \ref{thrm:main-symp-cocycle} is to check that the symplectic 1-cocycle data on $\pi^*(\cT)$ determines tight contact structures on the three strongly concave boundary components.

Recall that a $n$-fold simple branched covering of $(Y,K)$ is defined by a map
\[\rho: \pi_1(Y \setminus K) \rightarrow S_n\]
where each meridian in $\pi_1(Y \setminus K)$ is mapped to a transposition in $S_n$.  The kernel of $\rho$ corresponds to a covering space $\widehat{Y} \setminus \widehat{K}$ of $Y \setminus K$.  If $K$ is transverse to a contact structure $(Y,\xi)$, then (after a perturbation along the ramification locus $\widehat{K} \subset \widehat{Y}$) the contact structure pulls back to a contact structure $(\widehat{Y},\widehat{\xi})$ that is unique up to ambient isotopy.  A local model is given in \cite[Theorem 7.5.4]{Geiges}, for example.  

\begin{example}
Let $H$ denote the Hopf link and $T$ the right-handed trefoil.
\begin{enumerate}
    \item The link group of $H$ is $\pi_1(S^3 \setminus T) = \langle x,y | xy = yx \rangle = \pi_1(T^2)$.  Therefore, any homormorphism $\rho$ must send $x,y$ to commuting transpositions
    \item The knot group of $T$ is $\pi_1(S^3 \setminus T) \cong \langle x,y | xyx = yxy \rangle$.  The homomorphsim $\rho$ to $S_3$ defined by $\rho(x) = (12)$ and $\rho(y) = (23)$ defines an irregular 3-fold simple branched covering from $S^3$ to itself.
\end{enumerate}
\end{example}

\begin{proposition}
Let $(B^3,\xi_{std})$ denote the 3-ball with the standard tight contact structure and convex boundary.  Let $U$ denote a transverse unknot with $sl = -1$ and let $T$ denote a transverse, right-handed trefoil with $sl = 1$.
\begin{enumerate}
\item the 2-fold contact branched cover of $(B^3,\xi_{std})$ over $U$ is $(S^2 \times [0,1],\xi_{std})$, where $\xi_{std}$ is the unique tight contact structure.
\item the irregular 3-fold simple branched cover of $(B^3,\xi_{std})$ over $T$ is contactomorphic to $(B^3,\xi_{std})$;
\item Let $L \subset (S^3,\xi_{std})$ be a transverse link consists of the split union of maximal unknots; Hopf links (with both linking numbers) of maximal self-linking unknots; and maximal self-linking right-handed trefoils.  Let $(Y,\xi)$ be the contact branched covering over $L$ determined by a map $\rho: \pi_1(S^3 \setminus L) \rightarrow S_n$.  Then $(Y,\xi)$ is contactomorphic to $(S^3 \#_k S^1 \times S^2,\xi_{std})$ for some $k \geq 0$, where $\xi_{std}$ is the unique tight contact structure.
\end{enumerate}

\end{proposition}

\begin{proof}
Consider the unit ball $B$ in $\CC^2$.  This is a Stein domain and therefore its preimage under the branched covering $\pi: (z,w) \mapsto (z,w^2)$ is also Stein; in particular, the contact structure induced on the boundary is tight.  The ramification locus $\{w = 0\}$ intersects $S^3$ along a $sl = -1$ transverse unknot.  Take a $B^3$-neighborhood of this knot with convex boundary; the preimage under the branched covering is homeomorphic to $S^2 \times [0,1]$ and the restriction of the contact structure is therefore tight.  

Part (2) follows similarly using the cuspidal branched covering
\[\pi: (x,y) \mapsto (x, y^3 + xy)\]
whose branch locus in the codomain intersets $S^3$ along a right-handed trefoil.

Finally, note that we can choose $B^3$-neighborhoods of each split link component.  Then $Y$ is union of $n$ copies of a punctured $(S^3,\xi_{std})$ and components contactomorphic to $(S^2 \times [0,1],\xi_{std})$ or $(B^3,\xi_{std})$, glued along convex 2-spheres with tight neighborhoods.  This gluing preserves tightness, so since each component is tight, the global contact structure $(Y,\xi)$ is also tight.
\end{proof}

\subsection{First Chern class}

Let $\pi: (X,\omega) \rightarrow (\CP^2,\omega_{FS})$ be a symplectic branched covering as constructed in this paper.   We can interpret the first Chern class of $(X,\omega)$ geometrically in terms of the Euler class of the foliations $\cF_{\lambda} = \text{ker}(\betat_{\lambda})$ on $H_{\lambda}$ determined by a symplectic 1-cocycle $\{\betat_{\lambda}\}$.  The discussion follows \cite[Section 3]{LC-symp-cocycle}, where it was remarked that compact leaves minimize the (generalized) Thurston norm.

Let $B$ denote the zeros of $\beta_{\lambda}$ on $\Sigma = \partial H_{\lambda}$.  There is an intersection pairing
\[ \langle , \rangle: H_1(H_{\lambda},B) \times H_2(H_{\lambda}, \Sigma \setminus B) \rightarrow \ZZ \]
Choose a vector field $v_{\lambda}$ directing the kernel of $\beta_{\lambda}$ on $\Sigma$ and extend it to a section of $\text{ker}(\betat_{\lambda})$ over $H_{\lambda}$.  Generically, this section will vanish along a smooth 1-complex representing a well-defined Euler class $e(\cF_{\lambda}) \in H_1(H_{\lambda},B)$ of the foliation.  If $L$ is a compact leaf of the foliation $\cF_{\lambda}$, then evaluating $e(\cF_{\lambda})$ on $L$ gives its Euler characteristic.

\begin{proposition}
    The triple $(e(\cF_1),e(\cF_2),e(\cF_3))$ is a 1-cocyle that represents the first Chern class of $(X,\omega)$.
\end{proposition}

\begin{proof}
    To ensure that we obtain a 1-cocycle, we can assume that the vector fields satisfy the cocycle condition $v_1 + v_2 + v_3 = 0$.

    The first Chern class of $(X,\omega)$ is given by the formula
\[c_1(X,\omega) = \pi^*(c_1(\CP^2)) - PD(\cR)\]
where $\cR \subset X$ is the (smooth) ramification locus of $\pi$.  By construction, since the symplectic 1-cocycle is constructed via a branched covering over the standard trisection of $\CP^2$, every regular leaf $L$ of $\cF_{\lambda}$ is compact and is the $k$-fold branched cover of a disk leaf in $H_{\lambda} = S^1 \times D^2$.  Therefore, if the proposition holds with respect to the standard trisection of $\CP^2$, then it will hold in general by the Riemann-Hurwitz formula.

The first Chern class of $(\CP^2,\omega_{FS})$ has degree 3.  This can be represented by a triple $(C_1,C_2,C_3)$, where $C_{\lambda}$ is a core of the handlebody $H_{\lambda} \cong S^1 \times D^2$.  Since each compact leaf $L$ is a 2-disk, we see that $\langle C_{\lambda},L \rangle = \chi(L) = 1$.  
\end{proof}

\subsection{Realizing geometrically transverse surfaces}

\begin{proposition}
\label{prop:all-line-bundles}
    Let $(X,\omega)$ be a symplectic 4-manifold with rational symplectic form. The branch locus $\cR$ in the construction can be chosen such that every homology class $A \in H_2(X;\ZZ)$ with positive symplectic area is represented by a surface that is geometrically transverse to the trisection of $X$.
\end{proposition}

\begin{proof}
There exists a surjective map
\[\pi_1(\widetilde{H}_1,p_0) \rightarrow H_1(\widetilde{H}_1) \rightarrow H_2(X)\]
The first is given by the abelianization, while the second is due to the fact that $H_1(\widetilde{H}_1) \cong H_2(\widetilde{H}_1,\del \widetilde{H}_1)$ and there is a handle decomposition of $X$ such that compressing disks in $\widetilde{H}_1$ are precisely the cores of the 2-handles.

Let $A \in H_2(X;\ZZ)$ be a class of positive symplectic area and let $\zeta$ be a loop based at $p$ such that $[\zeta] \in \pi_1(\widetilde{H}_1,p)$ is a lift of $A$.  By assumption, $\zeta$ is nullhomologous in $Y_1,Y_3$ and therefore bounds slice surfaces in $Z_1,Z_3$ whose union is a surface representing $A$.  

As shown in the following subsection, there is a decomposition
\[ \zeta = l_1^n \mu_1 \rho_1 \cdots \mu_k \rho_k\]
where
\begin{enumerate}
    \item $l_1$ is loop in $\Ht_1$ transverse to the foliation determined by $\beta_1$
    \item each loop $\mu_i$ is tangent to the leaf of $\cF_1$ through the basepoint $p$,
    \item each $\rho_i$ is represented by a loop on $\Sigma$ that is homotopic in $\Ht_2$ to a loop in the leaf of $\cF_2$ through $p$ and that is contractible in $\Ht_3$.
\end{enumerate}

In a neighborhood of $p$, the defining 1-forms $\beta_1,\beta_2,\beta_3$ exact with local primitives $g_1,g_2,$ and $g_3 = -g_1 - g_2$.  Pick $2k$ points $p_1,\dots,p_{2k}$ is a small neighborhood of the basepoint $p$ such that
\begin{align*}
    g_1(p_i) > g_1(p_{i-1}) & &&\text{ for $i$ odd} \\
    g_2(p_i) > g_2(p_{i-1}) & && \text{ for $i$ even}\\
    g_3(p_i) > g_3(p_{i-1}) & &&\text{ for $i$ even}
\end{align*}

Note that since the points $\{p_i\}$ are in an arbitrary neighborhood of $p$, there is a unique identification between homotopy classes of loops based at $p$ and homotopy classes of paths from $p_j$ to $p_{j'}$ for any $j,j'$.

To build the surface, connect $p_{2i-1}$ to $p_{2i}$ by arcs of $\tau_2$ and $\tau_3$ homotopic to $\rho_i$ and connect $p_{2i}$ to $p_{2i+1}$ be an arc homotopic to $\mu_{i+1}$.  Finally, connect $p_{2k}$ to $p_1$ by a loop homotopic to $l_1^n \mu_1$.  If $n > 0$, then each arc of $\tau_{\lambda}$ can be assumed geometrically transverse to the foliation $\cF_{\lambda}$.  

From Stokes's Theorem, it follows that
\[\int_{\cK} \omega = \sum_{\lambda = 1}^3 \int_{\tau_{\lambda}} \beta_{\lambda}\]
Since all of the component arcs of $\tau_2,\tau_3$ are perturbations of flat paths, their contribution is $O(\epsilon)$, whereas the contribution of $\tau_1$ to his is approximately $n \int_{l_1} \beta_1$.  Therefore, the surface has positive symplectic area if and only if $n > 0$.
\end{proof}

\subsection{Fundamental group of $H_1 \smallsetminus \tau_1$}.  Since the branch locus $\cR$ is in bridge position, the tangle $(H_1,\tau_1)$ is trivial and the arcs of $\tau_1$ can be simultaneously isotoped into the boundary $\del H_1 = T^2$.  Consequently, the 3-manifold $H_1 \smallsetminus \nu(\tau_1)$ is a genus $b + 1$ handlebody, where $b$ is the bridge index, and
\[\pi_1(H_1 \smallsetminus \tau_1,p) = \langle l,\mu_1,\dots,\mu_b\rangle \]
where $l$ is the longitude generating $\pi_1(H_1,p)$ and the $\mu$'s are meridians of the arcs of $\tau_1$.  However, we require a different presentation of the fundamental group adapted to the foliations $\cF_{\lambda} = \text{ker}(\beta_{\lambda})$.

\begin{definition}
Let $Y$ be a compact 3-manifold with (smooth) foliation $\cF$.  We say that a path is {\it flat} if it can be homotoped, rel boundary, to lie in a leaf of $\cF$.  
\end{definition}

Consider the cover $\overline{H}_1 = D^2 \times \RR$ associated to the longitude $l$ and let $p_0$ denote the lift of $p$ with $\RR$-coordinate 0.  The tangle $(\overline{H}_1,\overline{\tau}_1)$ consists of an infinite number of arcs and $\pi_1(\overline{H_1} \smallsetminus \overline{\tau}_1,p_0)$ is freely generated by meridians of these arcs.  Define the following subgroups:
\begin{align*}
M_0 \coloneqq & \text{ the subgroup generated by meridians that can be represented by loops } \\ & \text{ in the disk $D \times \{0\}$.} \\
M_r \coloneqq & \text{ the subgroup generated by meridians that can be represented by loops } \\ & \text{ in the cylinder $D \times [0,r]$ or $D \times [r,0]$.} \\
M \coloneqq & \text{ the subgroup generated by meridians.}
\end{align*}
Note that each of these groups projects to a subgroup of $\pi_1(H_1 \smallsetminus \tau_1,p)$ and by abuse of notation we let it denote both groups.  Note that if $r > 0$, then $lM_rl^{-1}$ is contained in $M_{r +1}$, and if $r < 0$ then $l^{-1}M_r l$ is contained in $M_{r - 1}$.

By construction, we have a torus diagram for the branch locus $\cR$ such that $\tau_1$ projects onto $T^2$ as a collection $\cA$ of disjoint arcs.  This diagram also pulls back to a diagram on $S^1 \times \RR$ for $\overline{\tau}_1$.

\begin{lemma}
The subgroup $M_r$ is generated by a collection of meridians, one for each embedded arc in the torus diagram in the region $S^1 \times [0,r]$ or $S^1 \times [r,0]$
\end{lemma}

\begin{proof}
Immediate application of the Wirtinger presentation, where the tangle is inside the solid torus and we view it from the outside.
\end{proof}

\begin{lemma}
If there are no bridge points in the cylinder $S^1 \times [s,r]$, then $M_s = M_r$.
\end{lemma}

\begin{proof}
All of the arcs are embedded and end at bridge points.  Thus the number of arcs can only increase at a bridge point.
\end{proof}

\begin{proposition}
\label{prop:flat-decomposition}
There exist a subgroup $R \subset \pi_1(\Sigma \smallsetminus \overline{b},p)$ such that
\begin{enumerate}
\item $M_r \subset M_s \ast R$ if $0 < s < r$ or $r < s < 0$
\item $M \subset M_0 \ast R$
\item $R$ is invariant under conjugation by the longitude $l$
\item every loop in $R$ contracts in $H_3 \smallsetminus \tau_3$
\item every loop in $R$ maps to a flat loop in $H_2 \smallsetminus \tau_2$
\end{enumerate}
\end{proposition}

\begin{figure}[h!]
\centering
	\begin{subfigure}[b]{0.5\textwidth}
		\includegraphics[width=\textwidth]{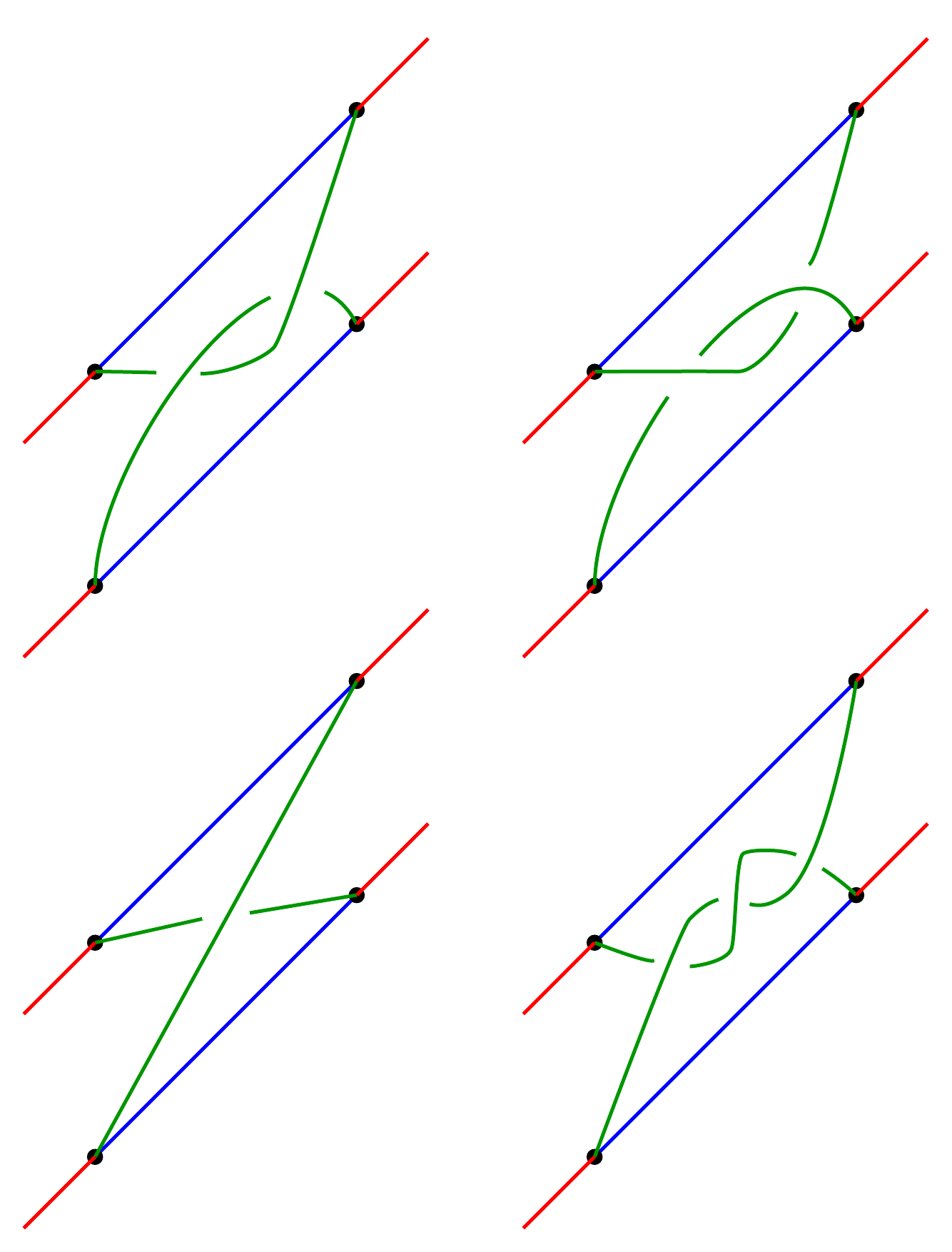}
		\caption{Four local models for $k = 1,3,2,-2$ (counterclockwise, starting in bottom left)}
		\label{fig:four-local-twists}
	\end{subfigure}
	\hfill
	\begin{subfigure}[b]{0.3\textwidth}
			\labellist
			\large\hair 2pt
			\pinlabel $a$ at 105 135
			\pinlabel $b$ at 110 250
			\pinlabel $d$ at 225 270
			\pinlabel $c$ at 240 160
		\endlabellist
		\includegraphics[width=\textwidth]{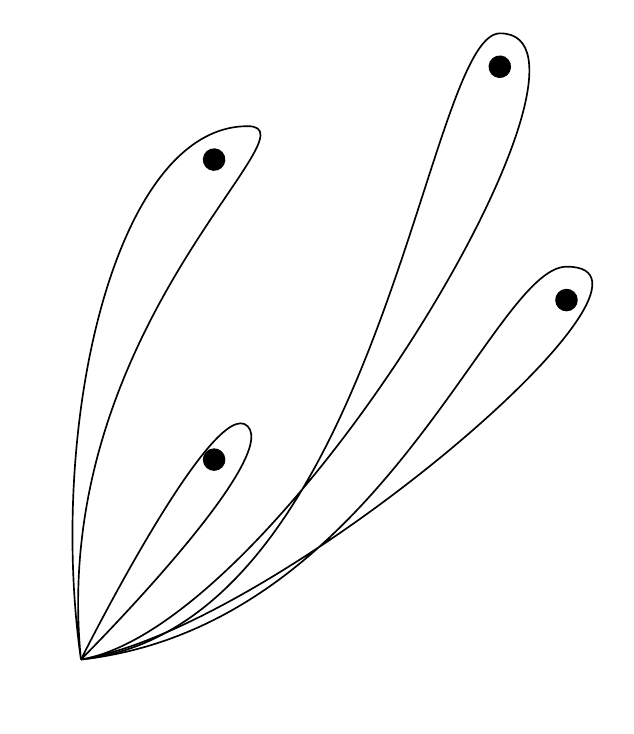}
		\caption{Paths in $\Sigmat$ that descend to meridians of the tangles $\tau_1,\tau_2,\tau_3$.  All oriented counterclockwise.}
		\label{fig:four-local-paths}
	\end{subfigure}
\caption{Local fundamental group calculation}
\label{fig:pi1-relations}
\end{figure}

\begin{proof}
We will prove part (1) by induction; then (2) clearly follows from (1). The properties of $R$ will come out of the construction.

For simplicity, we will assume that $0 < s < r$; the negative case is essentially identical. The group $M_r$ changes at a discrete collection of values of $r$, corresponding to two cases: (a) a crossing in the torus diagram, or (b) a local model of a half-twist.

{\bf Case (a): crossings}. Suppose that the crossing happens at level $r$. In the local picture \ref{fig:mini-stabilization-crossing}, we have three $\tau_1$-arcs and two bridge points.  We can choose loops $a,b$ in $\Sigma$ that descend to the Wirtinger meridians of the two arcs ending at the basepoints, so that $M_{r + \epsilon}$ is obtained from $M_{r - \epsilon}$ by adding the generator $b$.  Note that $a,b$ map to the same meridians of $\tau_{2}$ in $H_2$ and $\tau_3$ in $H_3$; thus $ab^{-1}$ contracts in $H_2 \smallsetminus \tau_2$ and in $H_3 \smallsetminus \tau_3$.  Include $\rho = ab^{-1}$ in $R$; clearly it satisfies properties (3)-(5) and $b = a \rho^{-1}$.

{\bf Case (b): half-twist}. In the local model, we have four bridge points. Choose four loops $a,b,c,d$ in $\Sigma$ based at $p$ as in Figure \ref{fig:four-local-paths}; we can assume they agree outside this local picture.  Pushing these loops into $H_{\lambda}$, they become meridians of $\tau_{\lambda}$. In $H_1$, they become the Wirtinger meridians of the four incident $\tau_1$ arcs.

First, we claim that each of the four loops maps to a flat meridian in $H_2$. Note that $a = c$ and $b = d$ in $H_2$. Then each $\cB$-arc in the projection intersects $\alpha_2$ transversely in a point. We can then represent these by a path that travels from $p$ along $\alpha_2$, then makes a small meridian loop, and returns to $p$ along $\alpha_2$.

Secondly, we claim that the following loops, depending on the value of $k$, contract in $H_3 \smallsetminus \tau_3$. See Figure \ref{fig:four-local-twists}. These are precisely the loops corresponding to the arcs of $\tau_3$.

\[
\begin{tabular}{c | c | c}
$k$ & $\rho_1$ & $\rho_2$ \\
\hline
1 & $ad$ & $cb$ \\
2 & $a^{-1}dab$ &$ cb^{-1}ab$\\
-2 & $ab^{-1}a^{-1}cab$ & $b^{-1}dbaba^{-1}$ \\
3 & $a^{-1}dab^{-1}ab$ &  $cb^{-1}a^{-1}bab$
\end{tabular}
\]

In each case, we have that
\[c = \rho_1 c_0 \qquad d = \rho_2 d_0\]
where $c_0,d_0$ lie in $M_{s}$ and $\rho_1,\rho_2$ are loops that bound in $H_3 \smallsetminus \tau_3$ and are flat loops in $H_2$.

We now let $R$ consist of all $l$-conjugates of all the relations at each half-twist. Thus, it satisfies property (3) by definition. Since each relation vanishes in $H_3 \smallsetminus \tau_3$, so does every conjugate and therefore every element in $R$. Finally we can homotope the longitude to lie in the leaf of $H_2$ bounded by $\alpha_2$. This implies that all $l$-conjugates of flat loops in $H_2$ are also flat.
\end{proof}

\begin{proposition}
\label{prop:meridian-factor}
Every element $\gamma \in \pi_1(H_1 \smallsetminus \tau_1,p)$ can be factored as
\[\gamma =l^p \gamma_0\]
where $p$ is an integer and $\gamma_0 \in M_0 \ast R$.
\end{proposition}

\begin{proof}
Clearly, $\gamma$ can be factored into a product of the longitude and the meridians of $\tau_1$. For any meridian $x$, we have that
\[x l = l (l^{-1} x l) = l \overline{x}\]
where $\overline{x} = (l^{-1} x l)$ is a different meridian. Consequently, we can move all longitudes to the front of the factorization. We are then left with a product of meridians $\gamma_0$, which lies in $M_0 \ast R$ by Proposition \ref{prop:flat-decomposition}.
\end{proof}

Let $\pi_1(\Ht_1,\overline{p})$ denote the set of homotopy classes of paths in $\Ht_1$ connecting any two points $p_i,p_j$ in the preimage of $p$.  It consists precisely of all lifts of all paths in $\pi_1(H_1 \smallsetminus \tau_1,p)$.  Let $\Mt_0$ and $\Rt$ denote the set of all lifts of $M_0$ and $R$.

\begin{proposition}
\label{prop:upstairs-factor}
Every element $\gamma \in \pi_1(\Ht_1,p_1)$ can be factored in $\pi_1(\Ht_1,\overline{p})$ into the form
\[\gamma = l_1^p \gamma_0\]
where $\gamma_0$ is a sequence of paths in $\Mt_0$ and $\Rt$. and $l_1$ is the lift of $l$ based at $p_1$.
\end{proposition}

\begin{proof}
Every loop in the branched cover projects to a loop in $H_1 \smallsetminus \tau_1$ based at $p$. We can use Proposition \ref{prop:meridian-factor} to factor it. Then each longitude $l$ lifts to $l_1$, the lift based at $p_1$, and each element of $M_0,R$ lifts to a path in $\Mt_0$ or $\Rt$.
\end{proof}

\section{Holomorphic trisections}

\begin{definition}
    Let $(X,J)$ be an almost-complex 4-manifold.  A {\it holomorphic trisection} of $(X,J)$ is a trisection such that
    \begin{enumerate}
        \item $J$ is integrable in an open neighborhood of the spine of $\cT$,
        \item there exists a neighborhood $U_{\lambda}$ of $H_{\lambda}$ and a holmorphic function $g_{\lambda}$ on $U_{\lambda}$ such that $H_{\lambda} = \{Re(g_{\lambda}) = 0\}$ and $Z_{\lambda}$ is locally defined as $\{\text{max}(Re(g_1),-Re(g_2)) \leq 0\}$; and
        \item the triple $(g_1,g_2,g_3)$ satisfy the cocycle relation
        \[Re(g_1) + Re(g_2) + Re(g_3) = 0\]
        on $U_1 \cap U_2 \cap U_3$, which is a neighborhood of the central surface $\Sigma$.
    \end{enumerate}
\end{definition}

\begin{example}
\begin{enumerate}
    \item The standard trisection of $\CP^2$ induces a holomorphic trisection, where $g_{\lambda} = \log z_{\lambda}$ for homogeneous coordinates $[z_1:z_2:z_3]$ on $\CP^2$.
    \item Zupan \cite{Zupan} shows that the curves
    \[C_d = \{z_1^{d-1}z_2 + z_2^{d-1}z_3 + z_3^{d-1}z_1 = 0\} \]
    are in bridge position with respect to the standard trisection of $\CP^2$.  Therefore, the projective surfaces
    \[V_d \{z_4^d + z_1^{d-1}z_2 + z_2^{d-1}z_3 + z_3^{d-1}z_1 = 0 \} \subset \CP^3 \]
    admit holomorphic trisections, obtained by pulling back the standard trisection of $\CP^2$ by the projection $V_d \rightarrow \CP^2$.
\end{enumerate}
\end{example}

\begin{definition}
    A {\it K\"ahler structure} on a holomorphic trisection is a K\"ahler form
    \[\omega = i \del \dbar \phi\]
    where $\phi$ is a K\"ahler potential such that
    \begin{enumerate}
        \item $\phi$ is also a defining function for $\partial \nu(H_1 \cup H_2 \cup H_3)$ making each boundary component strictly {\it pseudoconcave} with respect to the outward-normal coorientation.
        \item the induced contact structures on the boundary components are tight
    \end{enumerate}
\end{definition}

\begin{example}
\label{ex:Kahler-potential-cp2}
    On $\CP^2$, consider the K\"ahler potential
    \[\phi_N = \log( \frac{1}{N} |\mathbf{z}|^2 + |\mathbf{z}|^{2N})\]
    for $N \gg 0$. As $N \rightarrow \infty$, the 0-level sets converge to the standard trisection of $\CP^2$.
\end{example}

By following the proof of Theorem \ref{thrm:main-symp-cocycle}, combined with Theorem \ref{thrm:holomorphic-ribbon} proved below, we obtain the following.

\begin{theorem}
\label{thrm:kahler-structure-general}
    Let $(X,\omega)$ be a closed symplectic 4-manifold with rational symplectic form.  There exists a holomorphic trisection $\cT$ on $X$ with a K\"ahler structure $\omega'$ cohomologous to $\omega$.
\end{theorem}

\begin{proof}
    By Theorem \ref{thrm:holomorphic-ribbon}, the branch locus can be assumed holomorphic in some neighborhood of the spine.  Then, for $N$ sufficiently large, the K\"ahler potential in Example \ref{ex:Kahler-potential-cp2} defines a K\"ahler form and the superlevel set $\phi_n^{-1}([0,\infty))$ is a neighborhood of the spine in which the branch locus is holomorphic.  Pulling this back gives a K\"ahler structure on the holomorphic trisection of the branched cover. 
\end{proof}

    Note that the form $\pi^*(\omega_{FS})$ is degenerate along the ramification locus in $X$, but this can be eliminated by a $C^{\infty}$-small perturbation of the K\"ahler potential $\pi^*(\phi)$.

\subsection{Holomorphic ribbon}

Let $K \subset \CC^n$ be a subset.  The {\it rational hull} of $K$ is the set
\[\widehat{K}_{\QQ} := \{z \in \CC^n : |f(z)| \leq \text{max}_{w \in K} |f(w)| \text{ for all meromorphic functions $f$ on $\CC$}\} \]
A compact set $K$ is {\it rationally convex} if $\widehat{K}_{\QQ} = K$.  

\begin{lemma}
    Let $K$ be a piecewise smooth 1-complex in $\CC^2$.  Then $K$ is rationally convex and admits a Stein neighborhood basis.  
\end{lemma}

\begin{proof}
    Let $z$ be a point in the complement of $K$.  Consider the `projection from a point' map $\CP^{2}$ $\setminus z \rightarrow \CP^1$, which sends each point $y \in \CP^{2}$ $\setminus z$ to the unique line through $z$ and $y$.  The image of $K$ is a piecewise-smooth 1-complex.  In particular, it has measure 0.  Let $L$ be the projective line corresponding to some point in the complement and let $L'$ be some $C^{\infty}$-close projective line that misses $z$.  If $f$ is the defining equation for $L'$, then $1/f$ is a rational  function that has a pole along $L'$ and if $L'$ is sufficiently close to $L$, then $1/|f(z)|$ is greater than the maximum modulus of $1/f$ along $K$.
\end{proof}

Consequently, we can find a strictly plurisubharmonic function defining a tubular neighborhood of $\cK$.  Within this neighborhood, we can find a holomorphic function whose intersection with the spine of the trisection is arbitrarily close to $\cK$. 

The following is a standard fact in complex analysis, which for example can be proven by H\"ormander's weighted $L^2$ spaces \cite{Hormander} (see also \cite[Lemma 28]{Donaldson}).

\begin{lemma}
\label{lemma:stein-approx}
    Let $\phi: \Omega \rightarrow \RR$ be a smooth plurisubharmonic function on a domain $\Omega \subset \CC^2$ and let $\Delta_r = \phi^{-1}((-\infty,r])$ denote a sublevel set of $\phi$.  For each pair $r < s$, there exists a constant $K = K(r,s)$ such that if $f$ is any smooth, complex-valued function on $\Delta_s$, there exists a holomorphic function $\widetilde{f}$ on the interior of $\Delta_r$ such that
    \[ || f - \widetilde{f} ||_{C^1(\Delta_r)} \leq K \left(|| \dbar f ||_{C^0(\Delta_s)} + || \nabla (\dbar f) ||_{C^0(\Delta_s)} \right)\]
\end{lemma}

With this fact in hand, we can approximate the spine of the surface by a holomorphic curve.

\begin{theorem}
\label{thrm:holomorphic-ribbon}
    Let $\cK$ be a surface that is geometrically transverse with respect to the standard trisection of $\CP^2$.  There exists a neighborhood $U$ of the spine $\cS$, an isotopic surface $\cK'$ whose restriction to $U$ is holomorphic, and such that $\cK' \cap \cS$ is $C^{1}$-close to $\cK \cap \cS$.
\end{theorem}

\begin{proof}
    By assumption, the surface $\cK$ is holomorphic in some neighborhood of the bridge points $\{b_j\} = \cK \pitchfork \Sigma$.  Let $U_j$ denote a fixed open neighborhood of $b_j$ on which $\cK$ is holomorphic and defined as the vanishing locus of $\widetilde{B}_j$.  Let $\tau_{\lambda,j}$ be an arc of $\cK \cap H_{\lambda}$.  We can choose a holomorphic bidisk neighbohood $\Delta_{\lambda,j}$ of $\tau_{\lambda,j}$ such that $\{\Delta_{\lambda,j}\}$ and $\{U_j\}$ give an open cover of $\cK \cap \cS$.  By shrinking the open sets as necessary, we can assume that all triple intersections among the open cover are empty.  Finally, we choose a strongly pseudoconvex neighborhood $V$ of $\cK$ contained within this open cover.

    Along $H_{\lambda}$, have coordinates $z_{\lambda}$ and $\theta_{\lambda+1} = \text{Im} \log(z_{\lambda+1})$.  There exists real-valued functions $f_{\lambda,j},g_{\lambda,j}$ such that $\tau_{\lambda,j}$ is cut out by
    \[z_{\lambda} = f_{\lambda,j}(\theta_{\lambda+1}) + i g_{\lambda,j}(\theta_{\lambda+1})\]
    We can take real analytic approximations $\widetilde{f},\widetilde{g}$ of these defining functions and complexify.  The vanishing locus of 
    \[z_{\lambda} = \widetilde{f}_{\lambda,j} + i \widetilde{g}_{\lambda,j}\]
    is a holomorphic curve in some neighborhood of $\tau_{\lambda,j}$.  In particular, $\widetilde{f},\widetilde{g}$ are defined on some smaller bidisk $\Delta'_{\lambda,j} \subset \Delta_{\lambda,j}$, but by the Oka-Weil theorem we can arbitrarily approximate $\widetilde{f},\widetilde{g}$ by holomorphic functions on $\Delta_{\lambda,j}$.  Let $\widetilde{F}_{\lambda,j} = \widetilde{f}_{\lambda,j} + i \widetilde{g}_{\lambda,j}$.  Suppose that $\tau_{\lambda,j}$ is incident to the bridge point $b_k$.  We can glue the defining functions $\{\widetilde{F}_{\lambda,j}\}$ and $\{\widetilde{B}_k\}$ into a defining function $G$ on $V$ as follows.  Choose a fixed partition of unity $\{\phi_{\lambda,j},\phi_k\}$ subordinate to the open cover and define
    \[G = \phi_{\lambda,j}\widetilde{F}_{\lambda,j} + \phi_{k} \widetilde{B}_k\]
    Since all triple overlaps in the open cover are empty, on the double overlaps $\Delta_{\lambda,j} \cap U_k$ we have
    \[\phi = \phi_{\lambda,j} = 1 - \phi_k \]
    Then
    \begin{align*}
        |\dbar G| &= |\dbar \phi| \cdot |\widetilde{F} - \widetilde{B}| \\
        |\nabla \dbar G| &\leq |\nabla (\dbar \phi)| \cdot |\widetilde{F} - \widetilde{B}| + |\dbar \phi| \cdot |\nabla (\widetilde{F} - \widetilde{B})|
    \end{align*}
    Therefore, as $\widetilde{F}_{\lambda,j} \rightarrow \widetilde{B}_k$ on each $\Delta_{\lambda,j} \cap U_k$, we use Lemma \ref{lemma:stein-approx} to obtain a holomorphic function $\widetilde{G}$ on $V' \subset V$ whose vanishing locus $C^1$-approximates $\cK \cap \cS$.  We can then perturb $\cK$ to agree with the vanishing locus of $\widetilde{G}$, which yields the required surface.
\end{proof}

\begin{corollary}
\label{cor:hol-ribbon}
    Let $(X,\cT)$ be a holomorphic trisection constructed as a branched cover of $(\CP^2,\cT_{std})$ and let $\cK$ be a surface that is geometrically transvese with respect to $\cT$.  Then $\cK$ can be $C^1$-approximated by a surface $\cK$ that is holomorphic in an open neighborhood of the spine of $\cT$.
\end{corollary}

\begin{proof}
    Let $\pi:X \rightarrow \CP^2$ be the branched covering, which is holomorphic in a neighborhood of the spines.  After a $C^{\infty}$-small perturbation, we can assume that $\pi(\cK)$ is disjoint from the branch locus in $\CP^2$ in some neighborhood of the spine.  The image is geometrically transverse to the standard trisection of $\CP^2$.  By the previous theorem, we can approximate it by a holomorphic ribbon, then pull this back to a holomorphic ribbon in $X$. 
\end{proof}

\section{Hodge Theory}

Since we can assume that $(X,\omega)$ admits a integrable complex structure over the 2-skeleton, we can investigate complex geometric objects.

\subsection{Picard group}

The first important object is the Picard group of holomorphic line bundles.

\begin{definition}
\label{def:Picard}
    Let $(X,J)$ be an almost-complex 4-manifold with holomorphic trisection $\cT$.  The {\it Picard group} of $(X,J,\cT)$ is the group of holomorphic line bundles over the spine of $\cT$, up to biholomorphic equivalence, that admit smooth trivializations along $\partial \cS$.  The group structure is given by tensor product.
\end{definition}

One immediate question regarding the Picard group concerns the Hodge structure induced on $H^2_{DR}(X)$.  The first Chern class gives a map
\[c_1: Pic(X,J,\cT) \rightarrow  H^2(X;\ZZ) \rightarrow H^2(X;\CC)\] 
For K\"ahler surfaces, the image is precisely the integral $(1,1)$-classes.  However, unlike the K\"ahler case, we can modify the construction and assume that every class in $H^2(X;\ZZ)$ can be represented by a holomorphic line bundle.

\begin{proof}[Proof of Theorem \ref{thrm:all-line-bundles}]

Start with some class $\alpha \in H^2(X;\ZZ)$ such that $\int_X \omega \wedge \alpha > 0$.  Then by Proposition \ref{prop:all-line-bundles}, the Poincare-dual of $\alpha$ can be represented by a surface that is geometrically transverse to the holomorphic trisection of $(X,\omega)$.  Applying Corollary \ref{cor:hol-ribbon}, this surface can be $C^1$-approximated by a holomorphic curve in some neighborhood $\nu(\cS)'$ of the spine.  This is an effective divisor, which immediately determines a holomorphic line bundle whose first Chern class is Poincare-dual to the divisor.  

To prove the theorem, note that $H^2(X;\ZZ)$ is finitely-generated and we can assume it is generated by cohomology classes with positive symplectic area.  Therefore, we can realize a generating set as holomorphic line bundles, then composing by tensor product realize a holomorphic line bundle with any given first Chern class.
\end{proof}

\subsection{Smooth Hodge decomposition}

Let $(X,\omega,J)$ be the symplectic structure constructed in the main theorem.  Then $\cS$ is a compact Riemannian 4-manifold with boundary.  Let $d$ be the exterior derivative on $k$-forms and $d^*$ is adjoint with respect to the metric.  Define
\[\cD = (d + d^*) \qquad \Delta = \cD^2 = d d^* + d^* d\]
Define the space of {\it harmonic forms}
\[\cH^k_{\Delta} = \{\alpha \in \Omega^k(X) : \Delta \alpha = 0\} \]
and {\it harmonic fields}
\[\cH^k_{\cD} = \{\alpha \in \Omega^k(X) : (d + d^*) \alpha = 0\}\]
Note that harmonic fields are $d$- and $d^*$-closed, but harmonic forms are not necessarily $d$- and $d^*$-closed on a compact manifold with boundary.  Moreover, in the case of nonempty boundary, the spaces of harmonic forms and fields are in general larger than DeRham cohomology groups of the manifold itself.

Since $X$ is almost-complex, we get splittings $\Omega^k(X) = \bigoplus_{p+q = k} \Omega^{p,q}(X)$ of $k$-forms into their $(p,q)$-components.  We also have operators
\[\dbar: \Omega^{p,q}(X) \rightarrow \Omega^{p,q+1}(X) \qquad \del: \Omega^{p,q}(X) \rightarrow \Omega^{p+1,q}(X) \]
and their metric adjoints $\dbar^*,\del^*$ and Laplacians
\[\Delta_{\dbar} = (\dbar + \dbar^*)^2 \qquad \Delta_{\del} = (\del + \del^*)^2\]

Let $U$ be a connected open set on which $J$ is integrable.  Then we have a fiberwise equality
\[\Delta = 2 \Delta_{\dbar} = 2 \Delta_{\del}\]
The latter two Laplacians clearly preserve $(p,q)$-bigrading.  Consequently, we get decompositions
\[\cH^k_{\Delta} = \bigoplus_{p+q = k} \cH^{p,q}_{\Delta} \]
over $U$.  Complex conjugation gives an isomorphisms
\[\cH^{p,q}_{\Delta} \cong \cH^{q,p}_{\Delta}\]
and the Hodge star gives isomorphisms
\[\cH^{p,q}_{\Delta} \cong \cH^{n-p,n-q}_{\Delta} \]
Wedging with the K\"ahler form $\omega$ induces maps
\[ L_{\Delta}: \cH^{p,q}_{\Delta} \rightarrow \cH^{p+1,q+1}_{\Delta} \]
which has a metric adjoint $\Lambda_{\Delta}$.  The K\"ahler identifies imply that $\oplus \cH^{p,q}_{\Delta}$ is a $\mathfrak{sl}(2,\CC)$-modules.

However, we do not necessarily get a decomposition of $\cH^k_{\cD}$ into $(p,q)$-summands.  We let $H^{p,q}_d(X)$ denote the vector space of $d$-closed $(p,q)$-forms modulo $d$-exact $(p,q)$-forms.  On this subspace, we can partial recover the classical Hodge-Riemann bilinear relations.

\begin{theorem}[Hodge-Riemann bilinear relations]
\label{thrm:Hodge-Riemann}

Let $X$ be a $n$-dimensional complex manifold with strictly pseudoconcave boundary.  Let $\rho$ be a strictly plurisubharmonic function defining the boundary.  Suppose that $X$ admits a K\"ahler form with K\"ahler potential $\rho$ near $\partial X$.  Then there is a nondegenerate pairing
\[H^{p,q}_d(X) \times H^{q,p}_d(X) \rightarrow \CC\]
for $p+q \leq n-1$.    
\end{theorem}

\begin{proof}
Since the K\"ahler form is compatible with the boundary, we will show that we can replace $\omega$ with some $\omega'$ that vanishes along $\partial X$.  Then the pairing $\langle \alpha, \beta \rangle$ is given by composing the wedge product
\[H^{p,q}_d \times H^{q,p}_d \rightarrow H^{p+q,p+q}_d \subset H^{2(p + q)}_{DR}(X)\]
with the pairing
\[H^{2(p + q)}_{DR}(X) \times H^{2(n - p -q)}_{DR}(X,\partial X) \rightarrow H^{2n}_{DR}(X,\partial X) \cong \CC\]
In particular,
\[\langle \alpha, \beta \rangle = \int_{X} \alpha \wedge \beta \wedge (\omega')^{n-p-q}\]
which is well-defined since $\omega' \in H^2_{DR}(X,\partial X)$.

Let $\rho$ be the defining function and let $f$ be a function such that
\[f'(0) = f''(0) = 0\]
and $f'(t) > 0$ and $f''(t) \geq 0$ for $t > 0$ and $f(t) = t$ for $t \geq  \delta > 0$.  Then
\[ \omega' = i \del \dbar f(\rho)\]
is a 2-form such that $\omega'|_{\partial X} = 0$ but $(\omega')^n > 0$ on the interior of $X$.
\end{proof}

\begin{corollary}
    Let $\alpha \in H^{p,q}_{\dbar}(X)$ be a global holomorphic $(p,0)$-form for $p  \leq n-2$.  Then $\alpha$ is $d$-closed.

    In particular, if $n \geq 2$ then $X$ admits no nonconstant holomorphic functions.
\end{corollary}

\begin{proof}
    We have that
    \[\int_X |d \alpha|^2 dvol = \int_X d \alpha \wedge d \overline{\alpha} \wedge \omega^{n - p - 1} = \int_{\partial X} \alpha \wedge d \overline{\alpha} \wedge \omega^{n-p-1} = 0\]
    when $n - p \geq 2$.  This implies that $|d \alpha| = 0$ pointwise, hence $d \alpha = 0$.
\end{proof}

\subsection{Dolbeault cohomology}

By construction, the spine $\nu(\cS)$ is a K\"ahler surface with pseudoconcave boundary.  Consequently, we have a Fr\"olicher spectral sequence
\[H^{\bullet,\bullet}_{\dbar}(\nu(\cS)) \Rightarrow H^{\bullet}_{d}(\nu(\cS))\]
and for each $0 \leq k \leq 2n$ obtain the bound
\[\sum_{p+q = k} \text{dim} H^{p,q}_{\dbar}(\nu(\cS)) \geq \text{dim}H^k_{d}(\nu(\cS))\]
which is weaker than the Hodge decomposition for a closed K\"ahler surface.

Let $E$ be a holomorphic vector bundle over $\nu(\cS)$.  One can then define the Dolbeault cohomology groups $H^{p,q}_{\dbar}(\nu(S),E)$.  Moreover, since the boundary of $\nu(\cS)$ is strictly pseudoconcave, we can apply the classical Andreotti-Grauert theorem.

\begin{theorem}[Andreotti-Grauert \cite{AG-finite}]
    Let $E$ be a holomorphic vector bundle over $\nu(\cS)$.  The rank of $H^{p,0}(\nu(\cS),E)$ is finite for all $p$.
\end{theorem}

However, $q = 1$ is the critical dimension and there is no finiteness result.  With suitable $\dbar$-Neumann boundary conditions along pseudoconcave boundary, Epstein proved finite-dimensionality for $\dbar$-harmonic sections of holomorphic bundles for all $q$ \cite{Epstein1,Epstein2,Epstein3}.

Given a complex surface $X = X_+ \cup X_-$ decomposed along a real hypersurface $Y$ such that $X_+$ has pseudoconvex boundary and $X_-$ has pseudoconcave boundary, Andreotti and Hill proved a Meyer-Vietoris sequence relating the Dolbeault cohomologies of $X,X_+,X_-$ and the Kohn-Rossi cohomology of the CR structure induced on $Y$ \cite{Andreotti-Hill}.  The failure of a comparable sequence here may shed light on non-K\"ahler symplectic structures.


\bibliographystyle{amsalpha}
\bibliography{References}

\end{document}